\documentclass[leqno,final]{siamltex}
\usepackage{graphicx} 
\usepackage{amsmath,amstext,amssymb,bm}
\usepackage{leftidx}

\usepackage{xcolor} 
\usepackage{soul} 
\usepackage{tikz}
\usetikzlibrary{shapes,arrows}
\usepackage{mathrsfs}
\usepackage{epstopdf}
\usepackage{color}
\usepackage{multirow}
\usepackage{tabularx}
\usepackage[shortlabels]{enumitem}

\numberwithin{equation}{section}
\newtheorem{remark}{Remark}[section]

\allowdisplaybreaks[4]

\newcommand{\cQ}{\mathcal{Q}}

\newcommand{\E}{\mathcal{E}}

\renewcommand{\div}{\mbox{\rm div\,}}

\newcommand{\cF}{\mathcal{F}}

\newcommand{\mH}{\mathbb{H}}
\newcommand{\mP}{\mathbb{P}}
\newcommand{\mE}{\mathbb{E}}
\newcommand{\mV}{\mathbb{V}}

\newcommand{\eps}{\epsilon}
\newcommand{\veps}{\varepsilon}
\newcommand{\Ome}{\Omega}
\newcommand{\p}{\partial}
\newcommand{\nab}{\nabla}
\newcommand{\vu}{{\bf u}}

\newcommand{\vH}{{\bf H}}

\newcommand{\vL}{{\bf L}}

\newcommand{\vv}{{\bf v}}
\newcommand{\ve}{{\bf e}_{\vu}}

\newcommand{\vw}{{\bf w}}
\newcommand{\vA}{{\bf A}}
\newcommand{\vG}{{\bf G}}

\newcommand{\pphi}{\pmb{\phi}}

\begin{document}
	
	\title{Full moment error estimates in strong norms for numerical  approximations of  stochastic Navier-Stokes equations with multiplicative noise, Part I:  time discretization\thanks{This work was partially supported by the NSF grants DMS-2309626 and DMS-2012414.} }
	\markboth{XIAOBING FENG AND LIET VO}{Time discretization for stochastic Navier-Stokes equations}

	\author{Xiaobing Feng\thanks{Department of Mathematics, The University of Tennessee, Knoxville, TN 37996, U.S.A. (xfeng@utk.edu).}
		\and
		Liet Vo\thanks{
			School of Mathematical and Statistical Sciences, The University of Texas Rio Grande Valley, Edinburg, TX 78539, U.S.A. (liet.vo@utrgv.edu). This author was partially supported by the NSF grant DMS-2530211.}
	}
	
	\maketitle
	
	\begin{abstract} 
		This paper focuses on deriving optimal-order full moment error estimates in strong norms for both velocity and pressure approximations in the Euler-Maruyama time discretization of the stochastic Navier-Stokes equations with multiplicative noise. Additionally, it introduces a novel approach and framework for the numerical analysis of nonlinear stochastic partial differential equations (SPDEs) with multiplicative noise in general. The main ideas of this approach include establishing exponential stability estimates for the SPDE solution, leveraging a discrete stochastic Gronwall inequality, and employing a bootstrap argument.
	\end{abstract}
	
	
	
	\section{Introduction}	\label{sec-1}
	This paper is the first installment in a sequel \cite{FV2024} which is devoted to numerical analysis of the stochastic Navier-Stokes equations with multiplicative noise  (see \eqref{eq1.1}) that has been an open issue for the stochastic Navier-Stokes equations and many other nonlinear stochastic partial differential equations (SPDEs) arising from many scientific and engineering applications although a lot of progress has been made in the recent years, see 	\cite{BBM14}--\cite{FV2024}, \cite{Qui_2022, LV2021}, and the references therein. 
	
	The stochastic Navier-Stokes equations can be thought of as a random force perturbation of the deterministic Navier-Stokes equations for modeling turbulent flows \cite{mcmullen2022navier}.  Specifically, the random force has the form $ \delta\vG(\vu) \frac{d W(t) }{dt}$, where $\{W(t); t\geq 0\}$ denotes a real-valued Wiener process and its (formal) time derivative represents a white noise.  $\vG$ is a given vector-valued function  (see Section \ref{sec-3} for its precise definition), $\vu=\vu(x,t,\omega)$ denotes the velocity field, and $\delta>0$ is a constant which controls the strength of noise.  When $\vG$ is independent of $\vu$, the noise is said to be {\em additive}, otherwise, it is called {\em multiplicative}.  The latter case is the focus of this paper, and we refer the reader to \cite{BM2021,FV2024} and references therein for detailed expositions on numerical analysis of the stochastic Navier-Stokes equations with additive noise. 
	
	Although adapting and extending the vastly available numerical methods for the deterministic Navier-Stokes equations to their stochastic counterpart is relatively straightforward (see \cite{BBM14}--\cite{FV2024}, \cite{Qui_2022, LV2021}), however, their error analyses could not be imported nor adapted easily because the numerical analyses for the deterministic and stochastic PDEs are fundamentally different, this is especially true for nonlinear SPDEs with multiplicative noise. Three primary reasons cause additional difficulties with nonlinear SPDEs, including the prototypical stochastic Navier-Stokes equations. Firstly,  unlike in the deterministic case, in which one only needs to derive error estimates at a single sample,  the focus of the numerical analysis for SPDEs is to obtain error estimates for the quantities of stochastic interest (such as the expectation, variance, higher moments, etc.).  Secondly, the interplay between the nonlinearity  (in both drift and diffusion terms) and the stochasticity, which is often problem-dependent, is difficult to control. Finally, the lack of new numerical analysis techniques for SPDEs is another big hurdle because many well-known numerical analysis techniques for deterministic PDEs do not apply or simply do not work. All three types of difficulties arise strongly in the case of the stochastic Navier-Stokes equations with multiplicative noise (see \eqref{eq1.1}) as explained below.  
	
	Notice that the only nonlinearity in the (deterministic) Navier-Stokes equations comes from the convective term $\vu\cdot\nabla\vu$, which gives rise to the term $(\ve\cdot\nabla\vu, \ve)$ in the error equations, where $\ve$ denotes the error function for a velocity approximation. In the deterministic case, this term can be easily controlled by using Young's inequality and the Sobolev embedding (combined with the Gronwall's inequality approach). However, in the stochastic case, one needs to provide an error estimate for a quantity of stochastic interest. The simplest and most popular such quantity is the expectation. In that case,  $\vu$ and $\ve$ are random fields, and one must control the Lebesgue integral  $\mE[(\ve\cdot\nabla\vu, \ve)]:=\int_{\Omega} (\ve\cdot\nabla\vu, \ve)\, {\rm d}\mathbb{P}$. Due to its tri-factor structure, it is impossible to control this term without using higher (than second) moments of the error in strong norms.  To remedy and compromise,  instead of estimating the full moment $\mE[\cdot]$, one may consider estimating a partial moment 
	$\mE_{\Omega_m}[\cdot]:= \int_{\Omega_m} \cdot\, \, {\rm d}\mathbb{P}$  of  the error in strong norms for some 
	$\Omega_m \subset \Omega$ with $\Omega_m$ converging to $\Omega$ as $m\to \infty$ (see \cite{CP2012,BCP12}) for the choices of $\Omega_m$).  
	Indeed, all the current state-of-the-art numerical analysis results given in  \cite{BBM14}--\cite{CP2012} and \cite{Qui_2022}  were obtained using the above idea/approach, which was pioneered in \cite{BCP12,CP2012}. We refer the reader to a recent survey paper \cite{BPW2025} for more discussions about those techniques and results.  It should be noted that in \cite{bessaih2019strong, bessaih2022space}, the authors were able to derive an estimate for the complement  $\Omega\setminus \Omega_m$ of the sub-sample space $\Omega_m$ and consequently obtained a full second-moment error estimate for the velocity approximation under the assumption that the noise is almost additive. However, only a very weak logarithmic convergence rate was obtained in the case of multiplicative noise. 
	On the other hand, polynomial-order full moment error estimates in strong norms for the velocity approximations were obtained in \cite{BM2021,FV2024} in the case of additive noise, however, the analysis techniques used there strongly rely on the additive noise structure, hence, they can not be extended to the multiplicative-noise case due to the additional difficulties caused by the interplay between velocity field $\vu$ and noise $\frac{dW}{dt}$ in the random force term.  
	Lastly,  we felt that it is not possible (at least extremely difficult if possible) to derive polynomial-order full moment error estimates in strong norms for the stochastic Navier-Stokes equations using the improved deterministic numerical analysis approach/techniques pioneered in \cite{BCP12,CP2012} and refined in \cite{bessaih2019strong, bessaih2022space}.  Novel and new approaches/techniques must be developed to overcome the difficulties of deriving the desired polynomial-order full moment error estimates. 
	
	The primary goal of this and its companion paper \cite{FV2024b} is to address the above open issue and to provide a solution to it.  Specifically, we are able to establish optimal-order full moment error estimates in strong norms in the Euler-Maruyama time discretization of the stochastic Navier-Stokes equations here, while its fully discrete mixed finite element discretization will be considered in \cite{FV2024b}.  It should be noted that the error estimates are obtained not only for the velocity approximation but also for the pressure approximation.  Moreover, arbitrarily high-order moment and pathwise error estimates in the $L^2$-norm are also obtained.  Additionally,  we present a novel numerical analysis framework for nonlinear SPDEs with multiplicative noise that consists of four key components: (i) establishing exponential stability estimates for the SPDE solution; (ii) utilizing a discrete stochastic Gronwall's inequality; (iii) employing a bootstrap argument; 
	(iv) utilizing a special stochastic inf-sup condition/estimate.  We expect that this new framework is also applicable to other nonlinear stochastic PDEs (such as nonlinear stochastic Schr\"odinger equations), which will be considered in a forthcoming work, and the results will be presented elsewhere. 
	
	The remainder of this paper is organized as follows. In Section \ref{sec-2}, we introduce the function and space notation as well as a few useful facts to be used in this paper. In Section \ref{sec-3}, we present the stochastic Navier-Stokes equations, the structure assumptions on $\vG$, the variational (weak) solution concept and its properties, including the H\"older continuity in the energy norm and a crucial exponential stability estimate for the velocity field.  
	In Section \ref{sec-4}, we state the well-known Euler-Maruyama scheme for the stochastic Navier-Stokes equations and quote some stability estimates from \cite{BCP12}.  Finally, Section \ref{sec-5}, which is the main section, is devoted to the error analysis of the Euler-Maruyama scheme. The highlights of the section are proving optimal-order full moment error estimates in the energy norm for both velocity and pressure approximations and establishing arbitrarily high-order moment error estimates in the  $L^2$-norm for the velocity approximation.

	\section{Preliminaries and  useful facts}\label{sec-2}
	\subsection{Space notation}\label{sec-2.1}
	In this paper, standard function and space notation will be adopted. 
	Let $\vH^1_0(D)$ denote the subspace of $\vH^1(D)$ whose ${\mathbb R}^d$-valued functions have zero trace on $\p D$, and $(\cdot,\cdot):=(\cdot,\cdot)_D$ denote the standard $L^2$-inner product, with induced norm $\Vert \cdot \Vert$. We also denote ${\bf L}^p_{per}(D)$ and ${\bf H}^{k}_{per}(D)$ as the Lebesgue and Sobolev spaces of functions that are periodic with period $L$ in each coordinate direction for almost every $\mathbf{x} \in D$ { and have vanishing mean. } $C$ denotes a generic constant that is independent of the mesh parameters $h$ and $k$.
	
	Let $(\Omega,\cF, \{\cF_t\},\mP)$ be a filtered probability space with the probability measure $\mP$, the 
	$\sigma$-algebra $\cF$ and the continuous  filtration $\{\cF_t\} \subset \cF$. For a random variable $v$ 
	defined on $(\Omega,\cF, \{\cF_t\},\mP)$,
	${\mathbb E}[v]$ denotes the expected value of $v$. 
	For a vector space $X$ with norm $\|\cdot\|_{X}$,  and $1 \leq p < \infty$, we define the Bochner space
	$\bigl(L^p(\Omega;X); \|v\|_{L^p(\Omega;X)} \bigr)$, where
	$\|v\|_{L^p(\Omega;X)}:=\bigl({\mathbb E} [ \Vert v \Vert_X^p]\bigr)^{\frac1p}$.
	We also define 
	\begin{align*}
		{\mathbb H} := \bigl\{{\bf v}\in  \vL^2_{per}(D) ;\,\div {\bf v}=0 \mbox{ in }D\, \bigr\}\, , \quad 
		{\mathbb V} :=\bigl\{{\bf v}\in  \vH^1_{per}(D) ;\,\div {\bf v}=0 \mbox{ in }D \bigr\}\, .
	\end{align*}
	
	We recall from \cite{Girault_Raviart86} the (orthogonal) Helmholtz projection 
	${\bf P}_{{\mathbb H}}: \vL^2_{per}(D) \rightarrow {\mathbb H}$  and define the Stokes operator ${\bf A} := -{\bf P}_{\mathbb H} \Delta: {\mathbb V} \cap \vH^2_{per}(D) \rightarrow {\mathbb H}$. 
	
	\subsection{Some useful facts and inequalities} In this subsection, we cite some useful facts and inequalities that will be used in later sections. First,  we recall some well-known properties of the convective term of the Navier-Stokes equations (cf. \cite{Temam}).  It is easy to check that 
	\begin{align*}
		\bigl(\vu\cdot\nab\vv,\vv\bigr) =0  \qquad\forall \vu, \vv \in \vH^1_{per}(D).
	\end{align*}
	It follows from  H\"older's  and Ladyzhenskaya's inequalities  that
	\begin{align}\label{Ladyzhenskaya}
		\bigl|\bigl(\vu\cdot\nab\vv,\vw\bigr)\bigr| \leq  C_L\|\vu\|^{\frac12}_{\vL^2}\|\nab\vu\|^{\frac12}_{\vL^2}\|\nab\vv\|_{\vL^2} \|\vw\|^{\frac12}_{\vL^2}\|\nab\vw\|^{\frac12}_{\vL^2}
	\end{align}
	for some constant $C_L>0$. Moreover, we recall the following Sobolev inequality: 
	\begin{equation}\label{eqn2.2}
		\|\nab \vu\|_{\vL^2} \leq C_0\|\vA \vu\|_{\vL^2}
	\end{equation}
	for some constant $C_0>0$. 
	
	Next, we state the following property of the $\mathbb{R}$-valued Wiener process.  We refer the reader to  \cite{Ichikawa} for its high-dimensional generalization. 
	
	\begin{lemma}
		Let $\{W(t):t\geq 0\}$ be a $\mathbb{R}$-valued Wiener process, then there holds 
		\begin{align}\label{mean_wiener}
			\mE\Bigl[|W(t) - W(s)|^{2m}\Bigr] \leq C_m|t-s|^{m}\qquad\forall m \in \mathbb{N},
		\end{align}
		where $C_m$ is a positive constant and $C_1 =1$. 
	\end{lemma}
	
	Next, we quote the well-known Burkholder-Davis-Gundy inequality whose proof can be found in \cite[Theorem 2.4]{brzezniak1997stochastic}. 
	
	\begin{lemma} 
		For $1<r<\infty$, suppose $\phi \in L^r(\Ome;\mathcal{K})$. Then, there exists a constant $C_r >0$ such that
		\begin{align}\label{BDG}
			\mE\left[\sup_{0\leq t \leq T} \left\|\int_{0}^{t}\phi(s)\, dW(s)\right\|^r_{\mathcal{K}}\right] \leq C_r\mE\left[\left(\int_{0}^T \|\phi(s)\|^2_{\mathcal{K}}\, ds\right)^{\frac{r}{2}}\right],
		\end{align}
		where $\mathcal{K}$ is a separable Hilbert space.
	\end{lemma}
	
	Finally, we cite the following discrete stochastic Gronwall inequality from \cite[Theorem 1]{kruse2018discrete}, which will play a vital role in our error analysis.
	
	\begin{lemma}\label{Stochastic_Gronwall}
		Let $\left\{M_n\right\}_{n \in \mathbb{N}}$ be an $\left\{\mathcal{F}_n\right\}_{n \in \mathbb{N}}$-martingale satisfying $M_0=0$ on a filtered probability space $\left(\Omega, \mathcal{F},\left\{\mathcal{F}_n\right\}_{n \in \mathbb{N}}, \mathbb{P}\right)$. Let $\left\{X_n\right\}_{n \in \mathbb{N}},\left\{F_n\right\}_{n \in \mathbb{N}}$, and $\left\{G_n\right\}_{n \in \mathbb{N}}$ be sequences of nonnegative and adapted random variables with $\mathbb{E}\left[X_0\right]<\infty$ such that
		\begin{align}\label{ineq2.3}
			X_n \leq F_n+M_n+\sum_{k=0}^{n-1} G_k X_k \quad \text { for all } n \in \mathbb{N}.
		\end{align}
		Then, for any $q \in(0,1)$ and a pair of conjugate numbers $\alpha, \alpha' \in[1, \infty]$, i.e.,  $\frac{1}{\alpha}+\frac{1}{\alpha'}=1$, satisfying $q \alpha<1$,  there holds 
		\begin{align}\label{ineq2.4}
			\mathbb{E}\left[\sup _{0 \leq k \leq n} X_k^q\right] \leq\left(1+\frac{1}{1-\alpha q}\right)^{\frac{1}{\alpha}}\left\|\prod_{k=0}^{n-1}\left(1+G_k\right)^q\right\|_{L^{\alpha'}(\Omega)}\left(\mathbb{E}\left[\sup _{0 \leq k \leq n} F_k\right]\right)^q.
		\end{align}
	\end{lemma}
	
	\begin{remark}\label{rem2.1}
		\begin{enumerate}
			\item[(a)] Since the power $q\in (0,1)$,  inequality \eqref{ineq2.4} provides a sub-first-order moment bound for the random sequence $\{X_n\}$ although the order can be arbitrarily close to one. 
			\item[(b)] Let  $\{Y_n\}_{n\geq 0}$ be a nonnegative adapted sequence if the left-hand side of \eqref{ineq2.3}  is replaced by $X_n+Y_n$, it is easy to verify that inequality \eqref{ineq2.4} is still valid with $X_k + Y_k$ in the place of $X_k$.  This generalized version will also be used in Section \ref{sec-5}. 
		\end{enumerate}
	\end{remark}

	\section{Stochastic Navier-Stokes equations, weak solutions, and  their  properties}\label{sec-3}
	
	We consider  the  following time-dependent stochastic Navier-Stokes equations with multiplicative noise
	\begin{subequations}\label{eq1.1}
		\begin{alignat}{2} \label{eq1.1a}
			d\vu &=\bigl[ \nu\Delta \vu - (\vu\cdot \nab)\vu -\nabla p \bigr] dt +  \vG(\vu)d W(t)  &&\qquad\mbox{a.s. in}\, D_T,\\
			\div \vu &=0 &&\qquad\mbox{a.s. in}\, D_T,\label{eq1.1b}\\
			\vu(0)&= \vu_0 &&\qquad\mbox{a.s. in}\, D,\label{eq1.1d}
		\end{alignat}
	\end{subequations}
	where $D = (0, L)^2 \subset \mathbb{R}^2 \,$ represents a period/cell of the periodic domain in $\mathbb{R}^2$, $\vu$ and $p$ stand respectively for the velocity field and the pressure of the fluid, $\{W(t); t\geq 0\}$ denotes real-valued Wiener process. In addition, $\vG$ is a given vector-valued function (see below for its definition).  Here we seek periodic-in-space solutions $(\vu,p)$ with period $\L$, that is, 
	$\vu(t,{\bf x} + \L{\bf e}_i) = \vu(t,{\bf x})$ and $p(t,{\bf x}+\L{\bf e}_i)=p(t,{\bf x})$ 
	almost surely 
	and for almost every $(t, {\bf x})\in (0,T)\times \mathbb{R}^2$  and $i=1,2$, where 
	$\{\bf e_1,e_2\}$ denotes the canonical basis of $\mathbb{R}^2$.
	
	It should be noted that the reason for only considering the 2D case is due to a lack of theoretical PDE results for \eqref{eq1.1} in the 3D case, which is not surprising because its (simpler) deterministic counterpart is still a famous open (millennial) problem in the field of PDEs.  On the other hand,  we remark that the numerical scheme and some numerical results may be dimension-independent. 
	
	\subsection{Structure assumptions}\label{sec-3.1x}
	
	We introduce some structural assumptions on the nonlinear diffusion function $\vG$.  
	Let $\vG:   \vH^1_{per}(D)	\rightarrow \vH^1_{per}(D)$ such that 
	\smallskip
	
	\begin{enumerate}[(B1)]
		\item  There exists a constant $C_{\vG}>0$ such that $\|{\vG}(\vv) - \vG(\vw)\|_{\vL^2(D)} \leq C_{\vG}\|\vv-\vw\|_{\vL^2(D)}$ for any $\vv,\vw \in  \vH^1_{per}(D)$. 
		
		\item There exists a constant $K>0$ such that $\|\vG(\vv)\|_{\vL^2} \leq K$ for any $ \vv \in \vL^2_{per}(D)$.
		
		\item There exists a constant $L>0$ such that $\|\nab \vG(\vv)\|_{\vL^2} \leq  L\|\nab \vv\|_{\vL^2}$ for any $\vv \in    \vH^1_{per}(D).$     
	\end{enumerate}
	
	\smallskip
	\begin{remark}\label{rem3.1} 
		
		(a) 	(B1) requires that $\vG$ is Lipschitz in the $\vL^2_{per}(D)$-norm with the Lipschitz constant $C_\vG$.  Moreover, (B1) also implies that 
		$	\|\vG(\vv)\|_{\vL^2}\leq C_\vG \|\vv\|_{\vL^2}+ \|\vG(0)\|_{\vL^2}$ for any $\vv\in \vH^1_{per}(D). $ Hence, $\vG$ grows at most linearly in the $\vL^2_{per}(D)$-norm. 
		%
		
		(b) (B2) implies that $\vG(\vv)$ is uniformly bounded in $\vL^2_{per}(D)$  for any  $ \vv \in \vL^2_{per}(D)$.
		
		(c)  (B3) implies that $\vG$ is stable in the $\vH^1(D)$-norm. 
		
		(d)  Examples of $\vG$ that satisfy the assumptions (B1)--(B3) include  $\vG(\vv) = \sin(\vv)$, $\vG(\vv)=\cos(\vv)$, and $\vG(\vv)=\frac{|\vv|^2}{{1+|\vv|^2}}$.
	\end{remark}

	\subsection{Weak formulation and solutions}\label{sec-3.2x}
	In this subsection, we recall the concept of variational (weak) solution for \eqref{eq1.1} and refer the reader to \cite{Chow07,PZ1992,LRS03} for the detailed expositions, especially the proofs of its existence and uniqueness.
	
	\smallskip 
	
	\begin{definition}\label{def2.1} 
		Given $(\Omega,\cF, \{\cF_t\},\mP)$, let $W$ be an ${\mathbb R}$-valued Wiener process. 
		Suppose ${\bf u}_0\in L^2(\Omega; {\mathbb V})$.
		An $\{\cF_t\}$-adapted stochastic process  $\{{\bf u}(t) ; 0\leq t\leq T\}$ is called
		a variational solution of \eqref{eq1.1} if ${\bf u} \in  L^2\bigl(\Omega; C([0,T]; {\mathbb V})) 
		\cap L^2\bigl(\Ome; L^2((0,T);\vH^2_{per}(D) )\bigr)$, 
		and satisfies $\mP$-a.s.~for all $t\in (0,T]$
		\begin{align}\label{equu2.8a}
			\bigl({\bf u}(t),  {\bf v} \bigr) + \int_0^t  \nu\bigl(\nab {\bf u}(s),  &\nab {\bf v} \bigr) 
			\,  ds  + \int_0^t \big(\vu(s)\cdot\nab\vu(s),\vv\big)\, ds
			\\\nonumber
			&=({\bf u}_0, {\bf v}) + {  \Bigl(\int_0^t \vG(\vu(s))\, dW(s), {\bf v} \Bigr)}  \qquad\forall  \, {\bf v}\in {\mathbb V}\, . 
		\end{align}
	\end{definition}
	
	Definition \ref{def2.1} only defines the velocity $\mathbf{u}$ for \eqref{eq1.1}, 	its associated pressure $p$ is subtle to define (cf. \cite{LRS03}).  Moreover, it was shown in \cite{FPL2021,FV2024} that the pressure $p$ can be obtained as the distributional time-derivative of a regular function $P$ whose existence is guaranteed by the following theorem. 
	
	\begin{theorem}\label{thm 2.2}
		Let $\{{\bf u}(t) ; 0\leq t\leq T\}$ be the variational solution of \eqref{eq1.1}. There exists a unique adapted process 
		$P\in {L^2\bigl(\Omega; L^2(0,T; H^1_{per}(D))\bigr)}$ such that $(\mathbf{u}, P)$ satisfies 
		$\mP$-a.s.~for all $t\in (0,T]$
		\begin{subequations}\label{equu2.100}
			\begin{align}\label{equu2.10a}
				&\bigl({\bf u}(t),  {\bf v} \bigr) + \nu\int_0^t  \bigl(\nab {\bf u}(s), \nab {\bf v} \bigr) \, ds + \int_0^t \big(\vu(s)\cdot\nab\vu(s),\vv\big)\, ds
				- \bigl(  \div \mathbf{v}, P(t) \bigr)  \\
				&  \hskip 0.8in =({\bf u}_0, {\bf v}) 
				+  {\int_0^t  \bigl( {{\bf G}}(\vu(s)), {\bf v} \bigr)\, dW(s)}  \qquad \forall  {\bf v}\in \vH^1_{per}(D)\, , \nonumber \\ 
				& \bigl(\div {\bf u}, q \bigr) =0 \qquad\forall  q\in L^2_{per}(D) .  \label{equu2.10b}
			\end{align}
		\end{subequations}
	\end{theorem}
	
	We refer the interested reader to \cite[Theorem 1.3]{FPL2021} for a similar proof and to \cite{LRS03} for the derivation of  $p$ as the distributional time-derivative of $P$ (and $P$ as a time-average of $p$). 
	
	\subsection{Properties of weak solutions}\label{sec-3.2}
	
	In this subsection, we first quote some known stability estimates and H\"older continuity results for the variational solution. We then establish a new exponential stability estimate, which plays a crucial role in our error analysis later.
	
	The stability estimates quoted below were established in \cite{Breit, CP2012}. 
	
	\begin{lemma}\label{lemm3.3}
		Let $\vu$ be the solution defined in Definition \ref{def2.1}.  
		\begin{enumerate}[{\rm (a)}]
			\item Assume that $\vu_0 \in L^{\rho}\bigl(\Ome; \mV\bigr)$ for some $\rho \geq 2$. Then there holds
			\begin{align*}
				\mE\biggl[\Bigl(\sup_{0\leq t \leq T} \|\nab\vu(t)\|^2_{\vL^2} + \int_0^T \nu\|\nab^2\vu(t)\|^2_{\vL^2}\, dt\Bigr)^{\frac{\rho}{2}} \biggr] \leq C_{1,\rho},
			\end{align*}
			where $C_{1,\rho} = C(T,\rho)\mE\Bigl[\|\nab \vu_0\|^\rho_{\vL^2}\Bigr]$ is a positive constant.
			\item Assume that $\vu_0 \in L^{\rho}\bigl(\Ome;\mV\cap \vH^2(D) \bigr)\cap L^{5\rho}\bigl(\Ome,\mV\bigr)$ for some $\rho \geq 2$. Then we have
			\begin{align*}
				\mE\biggl[\Bigl(\sup_{0\leq t \leq T} \|\nab^2\vu(t)\|^2_{\vL^2} + \int_0^T \nu\|\nab^3\vu(t)\|^2_{\vL^2}\, dt\Bigr)^{\frac{\rho}{2}} \biggr] \leq C_{2,\rho},
			\end{align*}
			where $C_{2,\rho} = C(T, \rho) \mE\biggl[\Bigl(\|\vu_0\|^2_{\vH^2} + \|\nab\vu_0\|^{10}_{\vL^2}\Bigr)^{\frac{\rho}{2}}\biggr]$ is another positive constant.
		\end{enumerate}
	\end{lemma}
	
	Next, we cite the following high moment H\"older continuity estimate for the variational solution and refer the reader to \cite{Breit,CHP2012} for detailed proof. 
	
	\begin{lemma}\label{lemma2.2}
		Suppose ${\bf u}_0$ satisfies the assumptions in Lemma \ref{lemm3.3} (b) for some $r \geq 2$. Then there exists a constant $C \equiv C(D_T, \vu_0,r)>0$, such that the variational solution to problem \eqref{eq1.1} satisfies
		for $s,t \in [0,T]$
		\begin{align}\label{equu2.20a}
			\mE\Bigl[\|\vu(t) - \vu(s)\|^{r}_{\mV}\Bigr] \leq C_{2,r}|t-s|^{r\gamma}\qquad\forall \gamma \in \Bigl(0,\frac12\Bigr).
		\end{align}
	\end{lemma}
	
	It should be noted that ensuring the above H\"older continuity estimate is the only reason to confine our consideration to the case of periodic boundary conditions. In fact, all numerical results of this paper still hold for Dirichlet boundary conditions provided that the above H\"older continuity estimate holds for the SPDE solution. 
	
	\smallskip
	
	Next, we state and prove some exponential stability for the velocity solution $\vu$.
	
	\begin{lemma}\label{lemma_exp} 
		Let $\vu$ be the variational solution of \eqref{eq1.1} and  $\sigma_0 := \frac{1}{16K^2}$. Suppose $\vu_0\in L^2(\Ome;\mV)$ such that 	
		$\mE\left[\exp\left(4\sigma\|\vu_0\|^2_{\vL^2}\right)\right]<\infty$ for any $\sigma \in (0,\sigma_0]$ and $\vG$ satisfies Assumption (B2). Then there holds 
		\begin{align*}
			\sup_{0\leq t \leq T}		\mE\left[\exp\left({\sigma}\| \vu(t)\|^2_{\vL^2}\right)\right] + \mE\left[\exp\left(2\sigma\nu\int_0^T\|\nab\vu(s)\|^2_{\vL^2}\, ds\right)\right]&\leq C_3,
		\end{align*}
		where the constant $C_3=  \left(\frac12\mE\left[\exp\bigl(4\sigma \| \vu_0\|^2_{\vL^2} + 4\sigma K^2T\bigr)\right] +  1\right)e^{\frac{T}{2}}$.
	\end{lemma}
	
	\begin{proof}
		Let $\Phi(\vu(t)) = \|\vu(t)\|^2_{\vL^2}$, then
		\begin{align*}
			\Phi'(\vu) (\vv) := 2( \vu,  \vv),\qquad 
			\Phi''(\vu)(\vv,\vw) := 2( \vv,  \vw).
		\end{align*}
		Using It\^o's formula, we get
		\begin{align*}
			\| \vu(t)\|^2_{\vL^2} &= \| \vu_0\|^2_{\vL^2}  + 2\nu \int_0^t (\vu(\xi) ,\Delta \vu(\xi))\, d\xi \\ \nonumber
			&\qquad + 2 \int_0^t (\vu(\xi), \vG(\vu(\xi)))\, dW(\xi)  
			+ \int_0^t \| \vG(\vu(\xi))\|^2_{\vL^2}\, d\xi.
		\end{align*}
		Integrating by parts and using $(B2)$ yields
		\begin{align}\label{eq_2.511}
			&\| \vu(t)\|^2_{\vL^2} + 2\nu\int_0^t\|\nab\vu(s)\|^2_{\vL^2}\, ds \\\nonumber
			&\qquad 
			\leq \| \vu_0\|^2_{\vL^2} + 2\int_0^t(\vu(\xi),\vG(\vu(\xi)))\, dW(\xi) + K^2T.
		\end{align}
		For $\sigma>0$, let
		\begin{align*}
			Z_t := 4\sigma  \int_0^t (\vu(\xi), \vG(\vu(\xi)))\, dW(\xi).
		\end{align*}
		Then, $Z_t$ is a martingale for all $t \in [0,T]$ and its quadratic variation satisfies 
		\begin{align}\label{quadratic_variation11}
			\bigl<Z\bigr>_t &\leq 16\sigma^2 \int_0^t \|\vu(\xi)\|^2_{\vL^2}\|\vG(\vu(\xi))\|^2_{\vL^2}\, d\xi   \\ \nonumber
			&	\leq  16\sigma^2K^2 \int_0^t \|\vu(\xi)\|^2_{\vL^2}\, d\xi 
			\leq  \sigma \int_0^t \|\vu(\xi)\|^2_{\vL^2}\, d\xi,
		\end{align}
		where the last inequality above is obtained by using the assumption that $\sigma \leq \frac{1}{16K^2}$.
		
		Using the monotonicity of the exponential function, we get
		\begin{align*}
			&\exp\left({\sigma}\| \vu(t)\|^2_{\vL^2} + 2\sigma\nu\int_0^t\|\nab\vu(s)\|^2_{\vL^2}\, ds\right)\\ \nonumber
			&\qquad  \leq \exp\bigl(\sigma \| \vu_0\|^2_{\vL^2} + \sigma K^2T\bigr) \exp\left(\frac12\bigl[Z_t - \frac12\langle Z\rangle_t\bigr]\right) \exp\left(\frac14\langle Z\rangle_t \right).
		\end{align*}
		Using Young's inequality and \eqref{quadratic_variation11} we obtain
		\begin{align*}
			&\exp\left({\sigma}\| \vu(t)\|^2_{\vL^2} + 2\sigma\nu\int_0^t\|\nab\vu(s)\|^2_{\vL^2}\, ds\right)\\ \nonumber
			&\qquad  \leq \frac14\exp\bigl(4\sigma \| \vu_0\|^2_{\vL^2} + 4\sigma K^2T\bigr) +  \frac12\exp\left(Z_t - \frac12\langle Z\rangle_t\right) \\\nonumber
			&\quad\qquad+  \frac14\exp\left(\int_{0}^t \sigma\|\vu(\xi)\|^2_{\vL^2}\, d\xi \right),
		\end{align*}
		which implies that
		\begin{align*}
			&\exp\left({\sigma}\| \vu(t)\|^2_{\vL^2}\right) + \exp\left(2\sigma\nu\int_0^t\|\nab\vu(s)\|^2_{\vL^2}\, ds\right)\\ \nonumber
			& \qquad \leq \frac12\exp\bigl(4\sigma \| \vu_0\|^2_{\vL^2} + 4\sigma K^2T\bigr) +  \exp\left(Z_t - \frac12\langle Z\rangle_t\right) \\\nonumber
			&\quad \qquad+  \frac12\exp\biggl(\int_{0}^t \sigma\|\vu(\xi)\|^2_{\vL^2}\, d\xi \biggr).
		\end{align*}
		
		Using Jensen's inequality on the last term on the right-hand side of the above inequality,  we obtain
		\begin{align}\label{eq3.7}
			&\exp\left({\sigma}\| \vu(t)\|^2_{\vL^2}\right) + \exp\left(2\sigma\nu\int_0^t\|\nab\vu(s)\|^2_{\vL^2}\, ds\right)\\ \nonumber
			&\qquad  \leq \frac12\exp\bigl(4\sigma \| \vu_0\|^2_{\vL^2} + 4\sigma K^2T\bigr) +  \exp\left(Z_t - \frac12\langle Z\rangle_t\right) \\\nonumber
			&\quad \qquad+  \frac12\int_{0}^t \exp\left(\sigma\|\vu(\xi)\|^2_{\vL^2}\right)\, d\xi.
		\end{align}
		
		We note that the second term on the right side of \eqref{eq3.7} is an exponential martingale, so taking the expectation on \eqref{eq3.7}, we obtain
		\begin{align}
			&\mE\left[\exp\left({\sigma}\| \vu(t)\|^2_{\vL^2}\right)\right] + \left[\exp\left(2\sigma\nu\int_0^t\|\nab\vu(s)\|^2_{\vL^2}\, ds\right)\right]\\ \nonumber
			&\qquad  \leq \frac12\mE\left[\exp\bigl(4\sigma \| \vu_0\|^2_{\vL^2} + 4\sigma K^2T\bigr)\right] +  1 +  \frac12\int_{0}^t \mE\left[\exp\left(\sigma\|\vu(\xi)\|^2_{\vL^2}\right)\right]\, d\xi.
		\end{align}
		
		Applying the deterministic Gronwall inequality, we obtain
		\begin{align}\label{eq3.9}
			&\mE\left[\exp\left({\sigma}\| \vu(t)\|^2_{\vL^2}\right)\right] + \left[\exp\left(2\sigma\nu\int_0^t\|\nab\vu(s)\|^2_{\vL^2}\, ds\right)\right]\\ \nonumber
			&\qquad  \leq  \left(\frac12\mE\left[\exp\bigl(4\sigma \| \vu_0\|^2_{\vL^2} + 4\sigma K^2T\bigr)\right] +  1\right)e^{\frac{T}{2}}.
		\end{align}
		
		The proof is complete by taking $\sup_{0\leq t \leq T}$ to \eqref{eq3.9}.
	\end{proof}
	
	\begin{remark}
		The proof of Lemma \ref{lemma_exp} crucially uses Assumption (B2). A natural question is whether (B2) can be relaxed. The following trivial ODE example $dX(t)=X(t) dW(t)$ suggests the answer is negative because $\mE[\exp(|X(t)|^2)]=\infty$. 
	\end{remark}
	
	\section{Time approximation scheme}\label{sec-4} 
	In this section, we introduce the well-known (implicit) Euler-Maruyama time-stepping method for problem \eqref{eq1.1}, which is used as a prototypical time discretization to present our ideas for deriving full moment error estimates in strong norms. We remark that these ideas are also applicable to other time discretization methods, which will be presented in a forthcoming work. 
	
	\subsection{Definition of the scheme}
	Let $M>0$ be a positive integer and $k = \frac{T}{M}$ be the time step size. Define time points $t_j :=jk$ for $j=0,1,2,\cdots, M$.  Then,  $0=t_0<t_1 <t_2 <\cdots <t_M=T$ forms a uniform partition of the interval $[0,T]$.
	
	Let $\{\vu^n,p^n\}$ denote the approximate velocity and pressure solutions to \eqref{eq1.1} at $t_n$. The Euler-Maruyama time discretization for \eqref{eq1.1} seeks the pair $\{\vu^{n+1},p^{n+1}\}$ via the following algorithm:
	
	\smallskip
	\textbf{Algorithm 1.} 
	Let $\vu^0 = \vu_0$ be a given $\mV$-valued random variable. Find $\bigl(\vu^{n+1}, p^{n+1}\bigr) \in L^2(\Ome;\mV \times L^2_{per}(D))$ such that for any $\pphi \in \vH^1_{per}(D)$ and $\psi \in L^2_{per}(D)$, there hold $\mP$-a.s.
	\begin{subequations} \label{eqn4.1}
		\begin{align}\label{equu3.1}
			\bigl(\vu^{n+1} - \vu^{n},\pphi\bigr) + \nu k\,\bigl(\nab\vu^{n+1},\nab\pphi\bigr) &+ k\,\bigl(\vu^{n+1}\cdot\nab\vu^{n+1},\pphi\bigr)
			- k\, \bigl(p^{n+1},\div \pphi\bigr) \\\nonumber
			&= \bigl(\vG(\vu^n)\Delta W_{n+1},\pphi\bigr),\\
			\bigl(\div \vu^{n+1},\psi\bigr) &=0,
		\end{align}
	\end{subequations}
	where $\Delta W_{n+1} := W(t_{n+1})- W(t_n)$.
	
	To eliminate $p^{n+1}$ in the above system,  taking $\pphi \in \mV$ in \eqref{equu3.1}, we get
	\begin{align}\label{eq_reforms}
		\bigl(\vu^{n+1} - \vu^{n},\pphi\bigr) + \nu k\,\bigl(\nab\vu^{n+1},\nab\pphi\bigr) &+ k\,\bigl(\vu^{n+1}\cdot\nab\vu^{n+1},\pphi\bigr) 
		\\ \nonumber
		&	= \bigl(\vG(\vu^n)\Delta W_{n+1},\pphi\bigr), 
	\end{align}
	which is an equation for $\vu^{n+1}$ only. 
	
	Although it will not be addressed in this paper, we remark that spatial discretizations can be devised based on both \eqref{eqn4.1} and \eqref{eq_reforms}.  However,  they lead to very different fully discrete numerical methods for \eqref{eq1.1}, each has its own advantages and disadvantages (cf. \cite{Temam, CHP2012, FPL2021, FV2024,FV2024b}). 
	
	\subsection{Stability properties of the scheme}
	
	We quote the following stability estimates for $\{\vu^n\}$, which will be used in our error analysis later and whose proof can be found in \cite{BCP12}. 
	
	\begin{lemma}\label{stability_means}
		Let $\vu_0 \in L^{2^q}(\Ome;\mV)$ for an integer $1 \leq q <\infty$ be given, such that $\mE\bigl[\|\vu_0\|^{2^q}_{\mV}\bigr] \leq C$. Then there exists a constant $C_{4,q}= C(T, q, \vu_0)$ such that the following estimates hold:
		\begin{enumerate}[{\rm (i)}]
			\item $\displaystyle\mE\biggl[\max_{1\leq n \leq M}\|\vu^n\|^{2^q}_{\mV} + \nu k\sum_{n=1}^M \|\vu^n\|^{2^q-2}_{\mV}\|{\bf A}\vu^n\|^2_{\vL^2}\biggr] \leq C_{4,q}$.
			\item $\displaystyle\mE\Biggl[\biggl(\sum_{n=1}^M \|\vu^{n} - \vu^{n-1}\|^2_{\mV}\biggr)^q + \biggl(\nu k\sum_{n=1}^M \|{\bf A}\vu^{n}\|^2_{\mV}\biggr)^q\Biggr] \leq C_{4,q}$.
		\end{enumerate}
	\end{lemma}
	
	\section{Full moment error estimates in strong norms} \label{sec-5}
	This section is devoted to establishing the first main result of this paper, namely, to derive error estimates for the velocity approximation generated by Algorithm 1. These optimal estimates are obtained in full moment and strong norms. We note that partial moment error estimates in strong norms were obtained in \cite{BCP2013, CHP2012, CP2012} using quite different techniques from those of this paper, and the techniques of \cite{BCP2013,CHP2012,CP2012} were recently refined to obtain full moment error estimates of logarithmic orders in \cite{bessaih2019strong, bessaih2022space}. 
	
	\subsection{Basic (sub-second-order) full moment error estimates for the velocity approximation}\label{sec-5.1}
	We first derive our basic sub-second-moment error estimates in the energy norm for the velocity approximation.
	
	\begin{theorem}\label{theorem_semi_chapter5} 
		Let $\vu$ be the variational solution to \eqref{equu2.100}, $\{\vu^{n}\}_{n=1}^M$ be generated by  Algorithm 1, and $\sigma_0 = \frac{1}{16K^2}$. Assume $\vu_0 \in L^{4}\bigl(\Ome;\mV\cap \vH^2(D) \bigr)\cap L^{20}\bigl(\Ome,\mV\bigr)$ such that $\mE\left[\exp\left(4\sigma\|\vu_0\|^2_{\vL^2}\right)\right] < \infty$ for any $\sigma \in (0,\sigma_0]$, and $\vG$ satisfies Assumptions (B1)--(B3). Then, for any $0 < \gamma < \frac12$ and $0 < q \leq 1-\epsilon$ with some $0<\epsilon <<1$, there holds
		\begin{align}\label{equu310}
			&\left(\mE\left[\max_{1\leq n\leq M}\|\vu(t_n) - \vu^n\|^{2q}_{\vL^2}\right]\right)^{\frac{1}{2q}} \\\nonumber
			&\qquad\qquad\qquad+ \left(\mE\left[\left(k\sum_{n = 1}^M\|\nab(\vu(t_n) - \vu^n)\|^2_{\vL^2}\right)^q\right]\right)^{\frac{1}{2q}}
			\leq \tilde{C}_1\, k^{\gamma},
		\end{align}
		where $\tilde{C}_1 = C(\veps,q, \vu_0, T)$ is a positive constant.
	\end{theorem}
	
	\begin{proof}
		Since the proof is long and quite involved, we divide it into {\em seven} steps for presentation clarity.
		
		{\em Step 1:} We first derive the basic error equation. To the end, subtracting \eqref{eq_reforms} from \eqref{equu2.8a} and then choosing $\vv = \ve^{n+1}$ we obtain
		\begin{align}\label{error_eq}
			\nonumber	\bigl(\ve^{n+1}-\ve^n,\ve^{n+1}\bigr) + \nu k\|\nab\ve^{n+1}\|^2_{\vL^2} &= \nu \int_{t_n}^{t_{n+1}}\bigl(\nab(\vu(s)-\vu(t_{n+1})),\nab\ve^{n+1}\bigr)\, ds \\ 
			&+ \int_{t_n}^{t_{n+1}}\bigl(\vu^{n+1}\cdot\nab\vu^{n+1} - \vu(s)\cdot\nab\vu(s),\ve^{n+1}\bigr)\, ds\\  \nonumber
			&+\biggl(\int_{t_n}^{t_{n+1}}\bigl(\vG(\vu(s))-\vG(\vu^n))\bigr)\,dW(s),\ve^{n+1}\biggr)\\\nonumber
			&=: {\tt I + II + III}.
		\end{align}
		The first term on the left-hand side of \eqref{error_eq} can be easily controlled by using the formula $2(a,a-b) = \|a\|^2 - \|b\|^2 + \|a-b\|^2$. It remains to bound the right-hand side of \eqref{error_eq}, which will be done in the next two steps.
		
		\smallskip
		{\em Step 2:}  We first consider the terms ${\tt I}$ and ${\tt II}$.   By the Schwarz inequality, we obtain
		\begin{align}\label{equu311}
			{\tt I} \leq \nu\int_{t_n}^{t_{n+1}}\|\nab(\vu(t_{n+1}) - \vu(s))\|^2_{\vL^2}\, ds + \frac{\nu k}{4}\|\nab\ve^{n+1}\|^2_{\vL^2}.
		\end{align}
		Then, using the fact that $\bigl(\vu^{n+1}\cdot\nab\ve^{n+1},\ve^{n+1}\bigr) = 0$ due to \eqref{eqn2.2}, we can  rewrite the second term on the right-hand side of \eqref{error_eq} as 
		\begin{align}\label{equu312}
			{\tt II}		&=\int_{t_n}^{t_{n+1}}\bigl(\vu(t_{n+1})\cdot\nab\vu(t_{n+1}) - \vu(s)\cdot\nab \vu(s),\ve^{n+1}\bigr)\, ds \\\nonumber
			&\qquad+\int_{t_n}^{t_{n+1}}\bigl(\vu^{n+1}\cdot\nab\vu^{n+1} - \vu(t_{n+1})\cdot\nab \vu(t_{n+1}),\ve^{n+1}\bigr)\, ds \\\nonumber
			&=\int_{t_n}^{t_{n+1}}\bigl((\vu(t_{n+1}) - \vu(s))\cdot\nab\vu(t_{n+1}),\ve^{n+1}\bigr)\, ds\\\nonumber
			&\qquad+\int_{t_n}^{t_{n+1}} \bigl(\vu(s)\cdot\nab(\vu(t_{n+1}) - \vu(s)),\ve^{n+1}\bigr)\, ds\\\nonumber
			&\qquad- \int_{t_n}^{t_{n+1}}\bigl(\ve^{n+1}\cdot\nab\vu^{n+1},\ve^{n+1}\bigr)\, ds +\int_{t_n}^{t_{n+1}}\bigl(\vu(t_{n+1})\cdot\nab\ve^{n+1},\ve^{n+1}\bigr)\, ds \\\nonumber
			&=\int_{t_n}^{t_{n+1}}\bigl((\vu(t_{n+1}) - \vu(s))\cdot\nab\vu(t_{n+1}),\ve^{n+1}\bigr)\, ds\\\nonumber
			&\qquad+\int_{t_n}^{t_{n+1}} \bigl(\vu(s)\cdot\nab(\vu(t_{n+1}) - \vu(s)),\ve^{n+1}\bigr)\, ds\\\nonumber
			&\qquad- \int_{t_n}^{t_{n+1}}\bigl(\ve^{n+1}\cdot\nab\vu(t_{n+1}),\ve^{n+1}\bigr)\, ds \\\nonumber
			&=:{\tt II_1 + II_2 + II_3}.
		\end{align}
		
		Using \eqref{Ladyzhenskaya} and the Schwarz inequality, we get
		\begin{align*}
			{\tt II_1} 
			&\leq C\|\nab\vu(t_{n+1})\|^2_{\vL^2} \int_{t_n}^{t_{n+1}} \, \|\vu(s) - \vu(t_{n+1})\|^2_{\vH^1}\, ds + \frac{\nu k}{8}\|\nab\ve^{n+1}\|^2_{\vL^2},\\
			{\tt II_2} &\leq  C\int_{t_n}^{t_{n+1}} \, \|\vu(s) - \vu(t_{n+1})\|^2_{\vH^1}\|\nab\vu(s)\|^2_{\vL^2} \, ds + \frac{\nu k}{8}\|\nab\ve^{n+1}\|^2_{\vL^2}.
		\end{align*}
		
		To control ${\tt II_3}$, we add  the terms $\pm \int_{t_n}^{t_{n+1}}\nab\vu(s)\, ds$  and then use \eqref{Ladyzhenskaya} to get 
		\begin{align*}
			{\tt II_3} &= \int_{t_n}^{t_{n+1}}\bigl(\ve^{n+1}\cdot\nab(\vu(s)-\vu(t_{n+1})),\ve^{n+1}\bigr)\, ds - \int_{t_n}^{t_{n+1}}\bigl(\ve^{n+1}\cdot\nab\vu(s),\ve^{n+1}\bigr)\, ds\\\nonumber
			&\leq \|\ve^{n+1}\|^2_{\vL^4}\int_{t_n}^{t_{n+1}} \|\nab(\vu(s)-\vu(t_{n+1}))\|_{\vL^2}\, ds + \|\ve^{n+1}\|^2_{\vL^4}\int_{t_n}^{t_{n+1}} \|\nab\vu(s)\|_{\vL^2}\, ds\\\nonumber
			&\leq C\|\ve^{n+1}\|^2_{\vL^2} \int_{t_n}^{t_{n+1}}\|\nab(\vu(s) - \vu(t_{n+1}))\|^2_{\vL^2}\, ds \\\nonumber
			&\qquad+ \frac{4C_L^2}{\nu}\|\ve^{n+1}\|^2_{\vL^2}\int_{t_n}^{t_{n+1}} \|\nab\vu(s)\|^2_{\vL^2}\, ds + \frac{\nu k}{8} \|\nab\ve^{n+1}\|^2_{\vL^2}.
		\end{align*}
		Therefore, from the estimates of ${\tt II_1, II_2, II_3}$, we obtain
		\begin{align}\label{equu313}
			{\tt II} &\leq C\|\nab\vu(t_{n+1})\|^2_{\vL^2}\int_{t_n}^{t_{n+1}} \, \|\vu(s) - \vu(t_{n+1})\|^2_{\vH^1} \, ds  \\\nonumber
			&\qquad+C\int_{t_n}^{t_{n+1}} \, \|\vu(s) - \vu(t_{n+1})\|^2_{\vH^1}\|\nab\vu(s)\|^2_{\vL^2} \, ds \\\nonumber
			&\qquad + C\|\ve^{n+1}\|^2_{\vL^2}\int_{t_n}^{t_{n+1}}\|\nab(\vu(s) - \vu(t_{n+1}))\|^2_{\vL^2}\, ds\\\nonumber
			&\qquad +   \frac{4C_L^2}{\nu}\|\ve^{n+1}\|^2_{\vL^2}\int_{t_n}^{t_{n+1}} \|\nab\vu(s)\|^2_{\vL^2}\, ds + \frac{3\nu k}{8} \|\nab\ve^{n+1}\|^2_{\vL^2}.
		\end{align}
		
		\smallskip
		{\em Step 3:}  Next, we bound the noise term ${\tt III}$ as follows. 
		\begin{align}\label{equu314}
			{\tt III} 
			&= \biggl(\int_{t_n}^{t_{n+1}}\bigl(\vG(\vu(s))-\vG(\vu^n)\bigr)\,dW(s),\ve^{n+1} - \ve^n\biggr) \\\nonumber
			&\qquad+ \biggl(\int_{t_n}^{t_{n+1}}\bigl(\vG(\vu(s))-\vG(\vu^n)\bigr)\,dW(s),\ve^n\biggr)\\\nonumber
			&= \biggl(\int_{t_n}^{t_{n+1}}\bigl(\vG(\vu(s))-\vG(\vu(t_n))\bigr)\,dW(s),\ve^{n+1} - \ve^n\biggr) \\\nonumber
			&\qquad+\bigl((\vG(\vu(t_n)) - \vG(\vu^n))\Delta W_n, \ve^{n+1} - \ve^n\bigr) \\\nonumber
			&\qquad+ \biggl(\int_{t_n}^{t_{n+1}}\bigl(\vG(\vu(s))-\vG(\vu^n)\bigr)\,dW(s),\ve^n\biggr)\\\nonumber
			&\leq 2\biggl\|\int_{t_n}^{t_{n+1}}\bigl(\vG(\vu(s))-\vG(\vu(t_n))\bigr)\,dW(s)\biggr\|^2_{\vL^2} \\\nonumber
			&\qquad+ 2\|(\vG(\vu(t_n)) - \vG(\vu^n))\Delta W_{n+1}\|^2_{\vL^2} \\\nonumber
			&\qquad + \frac14\|\ve^{n+1} - \ve^n\|^2_{\vL^2} + \biggl(\int_{t_n}^{t_{n+1}}\bigl(\vG(\vu(s))-\vG(\vu^n)\bigr)\,dW(s),\ve^n\biggr)\\\nonumber
			&\leq 2\biggl\|\int_{t_n}^{t_{n+1}}\bigl(\vG(\vu(s))-\vG(\vu(t_{n}))\bigr)\,dW(s)\biggr\|^2_{\vL^2} + \frac14\|\ve^{n+1} - \ve^n\|^2_{\vL^2} \\\nonumber
			&\qquad + 2\|(\vG(\vu(t_n)) - \vG(\vu^n))\Delta W_{n+}\|^2_{\vL^2} \\\nonumber
			&\qquad+ \biggl(\int_{t_n}^{t_{n+1}}\bigl(\vG(\vu(s))-\vG(\vu^n)\bigr)\,dW(s),\ve^n\biggr).
		\end{align}
		
		\smallskip 
		{\em Step 4:}  We now can obtain the master inequality by combining the results of {\em Steps 1-3}. Indeed, 	substituting all the estimates for terms {\tt I, II, III} into \eqref{error_eq}, we obtain
		\begin{align}\label{eq_3.12}
			&\frac{1}{2}\Bigl[\|\ve^{n+1}\|^2_{\vL^2} - \|\ve^n\|^2_{\vL^2}\Bigr] +\frac14\|\ve^{n+1} - \ve^n\|^2_{\vL^2}+ \frac{3\nu k}{8}\|\nab \ve^{n+1}\|^2_{\vL^2}\\\nonumber
			&\quad \leq \nu \int_{t_n}^{t_{n+1}} \|\nab(\vu(t_{n+1}) - \vu(s))\|^2_{\vL^2}\, ds \\\nonumber
			&\qquad+ C\|\nab\vu(t_{n+1})\|^2_{\vL^2} \int_{t_n}^{t_{n+1}} \, \|\vu(s) - \vu(t_{n+1})\|^2_{\vH^1}\, ds  \\\nonumber&\qquad+C\int_{t_n}^{t_{n+1}} \, \|\vu(s) - \vu(t_{n+1})\|^2_{\vH^1}\|\nab\vu(s)\|^2_{\vL^2} \, ds \\\nonumber
			&\qquad + C\|\ve^{n+1}\|^2_{\vL^2}\int_{t_n}^{t_{n+1}}\|\nab(\vu(s) - \vu(t_{n+1}))\|^2_{\vL^2}\, ds\\\nonumber
			&\qquad +   \frac{4C_L^2}{\nu} \|\ve^{n+1}\|^2_{\vL^2}\int_{t_n}^{t_{n+1}}\|\nab\vu(s)\|^2_{\vL^2}\, ds  
			{+ 2k\|\vG(\vu(t_n)) - \vG(\vu^n)\|^2_{\vL^2} } \\\nonumber
			&\qquad+2\biggl\|\int_{t_n}^{t_{n+1}}\bigl(\vG(\vu(s))-\vG(\vu(t_{n}))\bigr)\,dW(s)\biggr\|^2_{\vL^2}  \\\nonumber
			&\qquad + 2\|(\vG(\vu(t_n)) - \vG(\vu^n))\Delta W_{n+1}\|^2_{\vL^2} 
			{ - 2k\|\vG(\vu(t_n)) - \vG(\vu^n)\|^2_{\vL^2} } \\\nonumber
			&\qquad + \biggl(\int_{t_n}^{t_{n+1}}\bigl(\vG(\vu(s))-\vG(\vu^n)\bigr)\,dW(s),\ve^n\biggr)\\ \nonumber
			&\quad \leq \nu \int_{t_n}^{t_{n+1}} \|\nab(\vu(t_{n+1}) - \vu(s))\|^2_{\vL^2}\, ds \\\nonumber
			&\qquad+ C\|\nab\vu(t_{n+1})\|^2_{\vL^2} \int_{t_n}^{t_{n+1}} \, \|\vu(s) - \vu(t_{n+1})\|^2_{\vH^1}\, ds  \\\nonumber&\qquad+C\int_{t_n}^{t_{n+1}} \, \|\vu(s) - \vu(t_{n+1})\|^2_{\vH^1}\|\nab\vu(s)\|^2_{\vL^2} \, ds \\\nonumber
			&\qquad + C\|\ve^{n+1}\|^2_{\vL^2} \int_{t_n}^{t_{n+1}}\|\nab(\vu(s) - \vu(t_{n+1}))\|^2_{\vL^2}\, ds\\\nonumber
			&\qquad +  C\|\ve^{n+1}-\ve^n\|^2_{\vL^2} \int_{t_n}^{t_{n+1}} \|\nab\vu(s)\|^2_{\vL^2}\, ds  \\\nonumber
			&\qquad +   \frac{4C_L^2}{\nu} \|\ve^{n}\|^2_{\vL^2} \int_{t_n}^{t_{n+1}}\|\nab\vu(s)\|^2_{\vL^2}\, ds + 2Ck\|\ve^n\|^2_{\vL^2} \\\nonumber
			&\qquad+2\biggl\|\int_{t_n}^{t_{n+1}}\bigl(\vG(\vu(s))-\vG(\vu(t_{n}))\bigr)\,dW(s)\biggr\|^2_{\vL^2}  \\\nonumber
			&\qquad + 2\|(\vG(\vu(t_n)) - \vG(\vu^n))\Delta W_{n+1}\|^2_{\vL^2} - 2k\|\vG(\vu(t_n)) - \vG(\vu^n)\|^2_{\vL^2} \\\nonumber
			&\qquad + \biggl(\int_{t_n}^{t_{n+1}}\bigl(\vG(\vu(s))-\vG(\vu^n)\bigr)\,dW(s),\ve^n\biggr),
		\end{align}
		{	where we added  the terms $\pm 2k\|\vG(\vu(t_n)) - \vG(\vu^n)\|^2_{\vL^2}$ on the right-hand side of the first inequality. This was done to get the correct setup for using the stochastic Gronwall inequality later.  Assumption (B1) was used to obtain the second inequality. }
		
		\smallskip
		{\em Step 5:} We  want to use the stochastic Gronwall inequality on \eqref{eq_3.12}. To the end, 
		we first apply the summation operator  $\sum_{n = 0}^{\ell}$ ($0\leq \ell <M$) to \eqref{eq_3.12} and use the fact that $\ve^0 = 0$ to get 
		\begin{align}\label{eq_3.13}
			\|\ve^{\ell+1}\|^2_{\vL^2} + \frac{3\nu k}{4} \sum_{n = 0}^{\ell} \|\nab \ve^{n+1}\|^2_{\vL^2} \leq F_{\ell} + M_{\ell} + \sum_{n = 0}^{\ell} \mathcal{G}_n \|\ve^n\|^2_{\vL^2},
		\end{align}
		where
		\begin{align*}
			F_{\ell} &:= \sum_{n = 0}^{\ell} \biggl[2\nu \int_{t_n}^{t_{n+1}} \|\nab(\vu(t_{n+1}) - \vu(s))\|^2_{\vL^2}\, ds \\\nonumber
			&\qquad+ C\|\nab \vu(t_{n+1})\|^2_{\vL^2} \int_{t_n}^{t_{n+1}} \|\vu(s) - \vu(t_{n+1})\|^2_{\vH^1}\, ds \\\nonumber
			&\qquad+ C \int_{t_n}^{t_{n+1}} \|\vu(s) - \vu(t_{n+1})\|^2_{\vH^1}\|\nab 
			\vu(s)\|^2_{\vL^2}\, ds \\\nonumber
			&\qquad + C\|\ve^{n+1}\|^2_{\vL^2} \int_{t_n}^{t_{n+1}} \|\nab(\vu(s) - \vu(t_{n+1}))\|^2_{\vL^2}\, ds \\\nonumber
			&\qquad + C\|\ve^{n+1}-\ve^n\|^2_{\vL^2} \int_{t_n}^{t_{n+1}} \|\nab\vu(s)\|^2_{\vL^2}\, ds \\\nonumber
			&\qquad+ 2\biggl\|\int_{t_n}^{t_{n+1}}\bigl(\vG(\vu(s))-\vG(\vu(t_n))\bigr)\,dW(s)\biggr\|^2_{\vL^2}\biggr],\\\nonumber
			M_{\ell} &:= \sum_{n = 0}^{\ell} Z_{n},\\\nonumber
			Z_{n} &:=2 \|\vG(\vu(t_n)) - \vG(\vu^n)\|^2_{\vL^2}|\Delta W_{n+1}|^2 
			{ - 2k \|\vG(\vu(t_n)) - \vG(\vu^n)\|^2_{\vL^2} }\\\nonumber
			&\qquad +\biggl(\int_{t_n}^{t_{n+1}}\bigl(\vG(\vu(s))-\vG(\vu^n)\bigr)\,dW(s),\ve^n\biggr),\\\nonumber
			\mathcal{G}_n &:= \frac{4C_L^2}{\nu}\int_{t_n}^{t_{n+1}} \|\nab\vu(s)\|^2_{\vL^2}\, ds  { + 2Ck.}
		\end{align*}
		
		Suppose that $\{M_{\ell}; \ell\geq 0\}$ is a martingale  (this fact will be verified later in {\em Step 7}), using the stochastic Gronwall inequality \eqref{ineq2.4} with $\alpha = 1+\eps, \beta = \frac{1+\eps}{\eps}$,  $0<q \leq 1 - \eps$, for some $\eps\in (0,1)$ (also see Remark \ref{rem2.1}(b)) to \eqref{eq_3.13}, we obtain
		\begin{align}\label{eq3.13}
			&\Bigl(\mE\bigl[\sup_{0 \leq \ell \leq M}\|\ve^{\ell}\|^{2q}_{\vL^2}\bigr]\Bigr)^{\frac{1}{2q}} + \left(\mE\left[\left(k\sum_{n = 0}^M\|\nab \ve^{n}\|^2_{\vL^2}\right)^q\right]\right)^{\frac{1}{2q}}\\\nonumber
			&\qquad \quad \leq \Bigl(1+ \frac{1}{1- \alpha q}\Bigr)^{\frac{1}{2\alpha q}}\biggl[\mE\Bigl[\exp\Bigl( \beta q\sum_{n = 0}^{M-1}\mathcal{G}_n\Bigr)\Bigr]\biggr]^{\frac{1}{2\beta q}} \,\Bigl(\mE\Bigl[\sup_{0 \leq \ell < M} F_{\ell}\Bigr]\Bigr)^{\frac12}.
		\end{align}
		
		\smallskip
		{\em Step 6:} The right-hand side of \eqref{eq3.13} involves two expectation factors.  We now show that the first factor is bounded by a constant and the second one contains a factor with the desired power of $k$. 
		
		{
			Next, we use Lemma \ref{lemma_exp} to bound the first factor as follows. First, we have
			\begin{align}\label{eq_3.14}
				\mE\left[\exp\left(\beta q\sum_{n = 0}^{M-1}\mathcal{G}_n\right)\right] &= \mE\left[\exp\left(\frac{\beta q 4C^2_L}{\nu}\int_{0}^T \|\nab\vu(s)\|^2_{\vL^2}\, ds\right)\right]e^{2C\beta qT }\\\nonumber
				&= \mE\left[\exp\left(2\sigma\nu\int_{0}^T \|\nab\vu(s)\|^2_{\vL^2}\, ds\right)\right]e^{2C\beta qT },
			\end{align}
			where $\sigma = \frac{2\beta q C_L^2}{\nu^2}>0$. 
			
			We choose $\epsilon \in (0,1)$ such that  $1>\eps > \frac{1}{2}\left(\sqrt{4+\frac{\nu^4}{32^2C_L^4K^4}}-\frac{\nu^2}{32C_L^2K^2}\right)$, which implies that $\sigma < \frac{1}{16K^2}$, hence,  $\sigma$ fulfills  the requirement of Lemma \ref{lemma_exp}. Thus,  from Lemma \ref{lemma_exp} we get 
			\begin{align*}
				\mE\left[\exp\left(\beta q\sum_{n = 0}^{M-1}\mathcal{G}_n\right)\right]  &\leq C_3e^{2C\beta qT }.
			\end{align*}
		}

		To bound the second factor, $\bigl(\mE\bigl[\sup_{0 \leq \ell < M} F_{\ell}\bigr]\bigr)^{\frac12}$, we first obtain
		\begin{align}\label{eqn5.11}
			\mE\left[\sup_{0 \leq \ell \leq M-1} F_{\ell}\right] &\leq C\sum_{n = 0}^{M-1}\biggl\{ \int_{t_n}^{t_{n+1}} \mE\bigl[\|\nab(\vu(t_{n+1}) - \vu(s))\|^2_{\vL^2}\bigr]\, ds\biggr\} \\ \nonumber 
			&\quad+ C\sum_{n = 0}^{M-1} \biggl\{ \int_{t_n}^{t_{n+1}} \mE\bigl[\|\vu(s) - \vu(t_{n+1})\|^2_{\vH^1}\|\nab \vu(t_{n+1})\|^2_{\vL^2}\bigr]\, ds \\\nonumber
			&\quad+ C \int_{t_n}^{t_{n+1}} \mE\bigl[\|\vu(s) - \vu(t_{n+1})\|^2_{\vH^1}\|\nab 
			\vu(s)\|^2_{\vL^2}\bigr]\, ds \biggr\}\\\nonumber
			&\quad + C\sum_{n = 0}^{M-1}\biggl\{\int_{t_n}^{t_{n+1}} \mE\bigl[\|\ve^{n+1}\|^2_{\vL^2}\|\nab(\vu(s) - \vu(t_{n+1}))\|^2_{\vL^2}\bigr]\, ds\biggr\} \\\nonumber
			&\quad + C\sum_{n = 0}^{M-1}\biggl\{\int_{t_n}^{t_{n+1}} \mE\bigl[\|\ve^{n+1}-\ve^n\|^2_{\vL^2}\|\nab\vu(s)\|^2_{\vL^2}\bigr]\, ds \biggr\} \\\nonumber
			&\quad+ 2\sum_{n = 0}^{M-1}\mE\biggl[\biggl\|\int_{t_n}^{t_{n+1}}\bigl(\vG(\vu(s))-\vG(\vu(t_n))\bigr)\,dW(s)\biggr\|^2_{\vL^2}\biggr]\\\nonumber
			&=:T_1 + T_2 +T_3+T_4+T_5.
		\end{align}
		It remains to bound $T_j$ for $j=1,2,3,4, 5$. 
		
		By Lemma \ref{lemma2.2}, we have
		\begin{align*}
			T_1&\leq C_{2,2}k^{2\gamma} \qquad\mbox{for any } 0<\gamma < \frac12.
		\end{align*}
		It follows from Lemmas \ref{lemma2.2} and \ref{lemm3.3}(a) that 
		\begin{align*}
			T_2 &\leq  C\sum_{n = 0}^{M-1} \biggl\{\int_{t_n}^{t_{n+1}}\left(\mE\left[\|\vu(s) - \vu(t_{n+1})\|^4_{\vH^1}\right]\right)^{\frac12}\left( \mE\bigl[\|\nab 
			\vu(t_{n+1})\|^4_{\vL^2}\bigr]\right)^{\frac12}\, ds \\\nonumber
			&\qquad+ C \int_{t_n}^{t_{n+1}} \left(\mE\left[\|\vu(s) - \vu(t_{n+1})\|^4_{\vH^1}\right]\right)^{\frac12} \left(\mE\bigl[\|\nab 
			\vu(s)\|^4_{\vL^2}\bigr]\right)^{\frac12}\, ds \biggr\}\\\nonumber
			&\leq  CC_{1,4}C_{2,4}k^{2\gamma}  \qquad\mbox{for any } 0<\gamma < \frac12.
		\end{align*}
		Similarly, by Lemmas \ref{lemma2.2},   \ref{lemm3.3}, and  \ref{stability_means}, we also can show
		\begin{align*}
			T_3 &\leq  C\left(C_{1,4} + C_{4,2}\right)C_{2,4}k^{2\gamma}  \qquad\mbox{for any } 0<\gamma < \frac12.
		\end{align*}
		
		To bound $T_4$, we use Lemmas \ref{lemma2.2},  \ref{lemm3.3}, and \ref{stability_means} (ii)  to get
		\begin{align*}
			T_4 &\leq Ck\mE\left[\sup_{s\in [0,T]}\|\nab \vu(s)\|^2_{\vL^2}\sum_{n = 0}^{M-1} \|\ve^{n+1} - \ve^n\|^2_{\vL^2} \right]\\\nonumber
			&\leq Ck\left(\mE\left[\left(\sup_{s\in [0,T]}\|\nab\vu(s)\|^2_{\vL^2}\right)^2\right]\right)^{\frac12}\left(\mE\left[\left(\sum_{n = 0}^{M-1}\|\ve^{n+1} - \ve^n\|^2_{\vL^2}\right)^2\right]\right)^{\frac12}\\\nonumber
			&\leq CC_{1,4}k\left(\mE\left[\left(\sum_{n = 0}^{M-1}\|\vu(t_{n+1}) - \vu(t_n)\|^2_{\vL^2}\right)^2 + \left(\sum_{n = 0}^{M-1}\|\vu^{n+1} - \vu^n\|^2_{\vL^2}\right)^2\right]\right)^{\frac12}\\\nonumber
			&\leq CC_{1,4}(C_{2,4} + C_{4,2})k.
		\end{align*}
		
		Finally, by the It\^o's isometry, Assumption (B1), and Lemma \ref{lemma2.2}, we get 
		\begin{align*}
			T_5 = 2\delta^2\sum_{n = 0}^{M-1}\mE\biggl[\int_{t_n}^{t_{n+1}}\|\vG(\vu(s))-\vG(\vu(t_n))\|^2_{\vL^2}\,ds\biggr] 
			\leq CC_{\vG}C_{2,2}k^{2\gamma}.
		\end{align*}
		Substituting the above estimates for $T_1, T_2, T_3,T_4,T_5$ into the right-hand side of \eqref{eqn5.11} and using \eqref{eq_3.14}, we obtain
		
		\begin{align}
			\Bigl(\mE\bigl[\sup_{0 \leq \ell \leq M}\|\ve^{\ell}\|^{2q}_{\vL^2}\bigr]\Bigr)^{\frac{1}{2q}} + \left(\mE\left[\left(k\sum_{n = 0}^M\|\nab \ve^{n}\|^2_{\vL^2}\right)^q\right]\right)^{\frac{1}{2q}} 
			\leq \tilde{C}_1 k^{\gamma},
		\end{align}
		where $\tilde{C}_1 = C_{\veps,q}C_3\left[C_{2,2} + C_{1,4}C_{2,4} + (C_{1,4} + C_{4,2})C_{2,4} + C_{1,4}(C_{2,4} + C_{4,2}) + C_{\vG}C_{2,2}\right]$.
		
		\smallskip
		{\em Step 7:} To complete the proof,  we still need to verify our claim in {\em Step 5} that  $\{M_{\ell}; \ell\geq 0\}$ is a martingale. To the end, we first use the It\^o's isometry,  the assumption $(B2)$, and the Burkholder-Davis-Gundy inequality \eqref{BDG} to get
		\begin{align}\label{eq_3.16}
			\mE[|M_{\ell}|] &\leq \sum_{n = 0}^{\ell} \mE[|Z_n|]\\\nonumber
			&\leq 4k\sum_{n=0}^{\ell} \mE\bigl[\|\vG(\vu(t_n)) - \vG(\vu^n)\|^2_{\vL^2}\bigr] \\\nonumber
			&\qquad+ \mE\left[\left|\sum_{n = 0}^{\ell} \left(\int_{t_{n}}^{t_{n+1}}\bigl(\vG(\vu(s)) - \vG(\vu^n)\bigr)\, dW(s), \ve^n\right)\right|\right]\\\nonumber
			&\leq C + \mE\left[\left|\sum_{n = 0}^{\ell} \left(\int_{t_{n}}^{t_{n+1}}\bigl(\vG(\vu(s)) - \vG(\vu^n)\bigr)\, dW(s), \ve^n\right)\right|\right]\\\nonumber
			&\leq C +  \left(\mE\left[\sum_{n = 0}^{M-1}\int_{t_{n}}^{t_{n+1}} \|\vG(\vu(s)) - \vG(\vu^n)\|^2_{\vL^2}\|\ve^n\|^2_{\vL^2}\, ds\right]\right)^{\frac12} \\\nonumber
			&<\infty.
		\end{align}
		
		In addition, for any $0\leq n \leq M-1$, using the martingale property of the It\^o integrals, we have
		\begin{align*}
			\mE[Z_n] &=  \mE\bigl[\|\vG(\vu(t_n)) - \vG(\vu^n)\|^2_{\vL^2}|\Delta W_{n+1}|^2 \bigr]  { -  k\mE\bigl[ \|\vG(\vu(t_n)) - \vG(\vu^n)\|^2_{\vL^2} \bigr] } \\\nonumber
			&\qquad \quad  + \mE\biggl[\biggl(\int_{t_n}^{t_{n+1}}\bigl(\vG(\vu(s))-\vG(\vu^n)\bigr)\,dW(s),\ve^n\biggr)\biggr]\\\nonumber
			&=k \mE \bigl[\|\vG(\vu(t_n)) - \vG(\vu^n)\|^2_{\vL^2} \bigr] 
			{ - k\mE \bigl[ \|\vG(\vu(t_n)) - \vG(\vu^n)\|^2_{\vL^2} \bigr]  } =0.
		\end{align*}
		Then,  the conditional expectation of $M_{\ell}$ given $\{M_n\}_{n=0}^{\ell-1}$ is
		\begin{align}\label{eq3.17}
			\mE\bigl[M_{\ell} | M_0, M_1, \cdots, M_n \bigr]  
			&= \mE\bigl[Z_0+Z_1 + \cdots + Z_n| M_0, M_1, \cdots, M_n \bigr] \\\nonumber
			&\qquad+ \mE \bigl[Z_{n+1} + \cdots+Z_{\ell}| M_0, M_1,\cdots, M_{n} \bigr]\\\nonumber
			&= M_{n} + \mE\bigl[Z_{n+1} + \cdots + Z_{\ell} \bigr] = M_n.
		\end{align}
		
		Thus,  we conclude that $\bigl\{M_{\ell}; \ell \geq 0\bigr\}$ is a martingale using \eqref{eq_3.16} and \eqref{eq3.17}. The proof is complete.
	\end{proof}

	\subsection{{High moment and pathwise} error estimates for the velocity approximation}
	The basic error estimate obtained in Theorem \ref{theorem_semi_chapter5} implies  a strong convergence in the $L_{\omega}^{r}L_t^{\infty}L_x^2$- and $L_{\omega}^{r}L_t^{2}H_x^1$-norm for $0<r<2$.  We note that these sub-quadratic moment estimates are consequences of using the stochastic Gronwall inequality \eqref{ineq2.4} (also see Remark \ref{rem2.1}(b)). However, we show below, which is the goal of this subsection,  that such a limitation can be overcome by a bootstrap argument to obtain arbitrarily high-order moment estimates, which in turn infers a strong convergence in the $L_{\omega}^{r}L_t^{\infty}L_x^2$-norm for all $0<r<\infty$.  
	
	\begin{theorem}\label{theorem_semi_higher_moment} 
		Let $\vu$ be the variational solution to \eqref{equu2.100}, $\{\vu^{n}\}_{n=1}^M$ be generated by  Algorithm 1, and $\sigma_0 = \frac{1}{16K^2}$. Assume   
		$\vu_0 \in L^{2^r}\bigl(\Ome;\mV\cap \vH^2(D) \bigr)\cap L^{2^r5}\bigl(\Ome,\mV\bigr)$ for $r\geq 1$ such that $\mE\left[\exp\left(4\sigma\|\vu_0\|^2_{\vL^2}\right)\right] < \infty$ for any $\sigma \in (0,\sigma_0]$, and $\vG$ satisfies Assumptions (B1)--(B3). Then, for any {real number} $0 < \gamma < \frac12$ and $0 < q \leq 1-\epsilon$ with some $0<\epsilon <<1$, there holds
		\begin{align}\label{eq3.16}
			\left(\mE\left[\max_{1\leq n\leq M}\|\vu(t_n) - \vu^n\|^{2^rq}_{\vL^2}\right]\right)^{\frac{1}{2^rq}} 
			\leq \tilde{C}_2\, k^{\gamma}
		\end{align}
		for some constant  $\tilde{C}_2 = C(\veps,q,r,\vu_0,T)>0$.
	\end{theorem}
	
	\begin{proof}  
		Again, since the proof is long, we present it in six steps.  Notice that the case $r=1$ was proved in Theorem \ref{theorem_semi_chapter5}.   We begin with proving  \eqref{eq3.16} for $r=2$ and then complete the proof by the induction argument. 
		
		\smallskip
		{\em Step 1:} 
		Multiplying the error inequality \eqref{eq_3.12} by $\|\ve^{n+1}\|^2_{\vL^2}$ and using the identity $2a(a-b) = a^2 - b^2 + (a-b)^2$, we obtain
		\begin{align}\label{eq_3.17}
			&\frac14\Bigl[\|\ve^{n+1}\|^4_{\vL^2} - \|\ve^n\|^4_{\vL^2}\Bigr] + \frac14\Bigl(\|\ve^{n+1}\|^2_{\vL^2} - \|\ve^n\|^2_{\vL^2}\Bigr)^2 \\ \nonumber
			&\hskip 0.75in + \frac14\|\ve^{n+1} - \ve^n\|^2_{\vL^2}\|\ve^{n+1}\|^2_{\vL^2} + \frac{3\nu}{8}\|\nab\ve^{n+1}\|^2_{\vL^2}\|\ve^{n+1}\|^2_{\vL^2}\\\nonumber
			&\qquad \leq \biggl[ \nu \int_{t_n}^{t_{n+1}} \|\nab(\vu(t_{n+1}) - \vu(s))\|^2_{\vL^2}\, ds \\\nonumber
			&\qquad \qquad  + C\int_{t_n}^{t_{n+1}} \, \|\vu(s) - \vu(t_{n+1})\|^2_{\vH^1}\|\nab\vu(t_{n+1})\|^2_{\vL^2} \, ds  \\\nonumber
			&\qquad\qquad +C\int_{t_n}^{t_{n+1}} \, \|\vu(s) - \vu(t_{n+1})\|^2_{\vH^1}\|\nab\vu(s)\|^2_{\vL^2} \, ds\biggr]\|\ve^{n+1}\|^2_{\vL^2} \\\nonumber
			&\qquad\qquad + C\int_{t_n}^{t_{n+1}}\|\ve^{n+1}\|^2_{\vL^2}\|\nab(\vu(s) - \vu(t_{n+1}))\|^2_{\vL^2}\, ds\|\ve^{n+1}\|^2_{\vL^2}\\\nonumber
			&\qquad \qquad +  \frac{4C_L^2}{\nu}\int_{t_n}^{t_{n+1}} \|\ve^{n+1}\|^4_{\vL^2}\|\nab\vu(s)\|^2_{\vL^2}\, ds  \\\nonumber
			&\qquad\qquad+2\biggl\|\int_{t_n}^{t_{n+1}}\bigl(\vG(\vu(s))-\vG(\vu(t_{n}))\bigr)\,dW(s)\biggr\|^2_{\vL^2}\|\ve^{n+1}\|^2_{\vL^2}  \\\nonumber
			&\qquad \qquad + 2\|(\vG(\vu(t_n)) - \vG(\vu^n))\Delta W_{n+1}\|^2_{\vL^2}\|\ve^{n+1}\|^2_{\vL^2}  \\\nonumber
			&\qquad \qquad + \biggl(\int_{t_n}^{t_{n+1}}\bigl(\vG(\vu(s))-\vG(\vu^n)\bigr)\,dW(s),\ve^n\biggr)\|\ve^{n+1}\|^2_{\vL^2}\\\nonumber
			&=: L_1 + L_2 + L_3 + L_4 + L_5 + L_6.
		\end{align}
		
		\smallskip
		{\em Step 2:} We want to control terms $L_j$ for $j=1,2,\cdots,6$.  To bound $L_1$,  rewriting $\|\ve^{n+1}\|^2_{\vL^2}= \bigl( \|\ve^{n+1}\|^2_{\vL^2}-\|\ve^n\|^2_{\vL^2} \bigr) +\|\ve^n\|^2_{\vL^2}$ and using the discrete Young's inequality, we get
		\begin{align*}
			L_1 &= \Bigl(\|\ve^{n+1}\|^2_{\vL^2}- \|\ve^n\|^2_{\vL^2}\Bigr) \biggl[ \nu \int_{t_n}^{t_{n+1}} \|\nab(\vu(t_{n+1}) - \vu(s))\|^2_{\vL^2}\, ds \\\nonumber
			&\qquad+ C\int_{t_n}^{t_{n+1}} \, \|\vu(s) - \vu(t_{n+1})\|^2_{\vH^1}\|\nab\vu(t_{n+1})\|^2_{\vL^2} \, ds  \\\nonumber
			&\qquad+C\int_{t_n}^{t_{n+1}} \, \|\vu(s) - \vu(t_{n+1})\|^2_{\vH^1}\|\nab\vu(s)\|^2_{\vL^2} \, ds\biggr] 
			\\\nonumber
			&\qquad +  \|\ve^{n}\|^2_{\vL^2}\biggl[ \nu \int_{t_n}^{t_{n+1}} \|\nab(\vu(t_{n+1}) - \vu(s))\|^2_{\vL^2}\, ds \\\nonumber
			&\qquad+ C\int_{t_n}^{t_{n+1}} \, \|\vu(s) - \vu(t_{n+1})\|^2_{\vH^1}\|\nab\vu(t_{n+1})\|^2_{\vL^2} \, ds  \\\nonumber
			&\qquad+C\int_{t_n}^{t_{n+1}} \, \|\vu(s) - \vu(t_{n+1})\|^2_{\vH^1}\|\nab\vu(s)\|^2_{\vL^2} \, ds\biggr] 
			\\\nonumber
			&\leq Ck \biggl[  \int_{t_n}^{t_{n+1}} \|\nab(\vu(t_{n+1}) - \vu(s))\|^2_{\vL^2}\, ds \\\nonumber
			&\qquad+ C\int_{t_n}^{t_{n+1}} \, \|\vu(s) - \vu(t_{n+1})\|^2_{\vH^1}\|\nab\vu(t_{n+1})\|^2_{\vL^2} \, ds  \\\nonumber
			&\qquad+C\int_{t_n}^{t_{n+1}} \, \|\vu(s) - \vu(t_{n+1})\|^2_{\vH^1}\|\nab\vu(s)\|^2_{\vL^2} \, ds\biggr] + \frac{1}{32}\Bigl(\|\ve^{n+1}\|^2_{\vL^2} - \|\ve^n\|^2_{\vL^2}\Bigr)^2 \\\nonumber
			&\qquad +  C\biggl[  \int_{t_n}^{t_{n+1}} \|\nab(\vu(t_{n+1}) - \vu(s))\|^4_{\vL^2}\, ds \\\nonumber
			&\qquad+ C\int_{t_n}^{t_{n+1}} \, \|\vu(s) - \vu(t_{n+1})\|^4_{\vH^1}\|\nab\vu(t_{n+1})\|^4_{\vL^2} \, ds  \\\nonumber
			&\qquad+C\int_{t_n}^{t_{n+1}} \, \|\vu(s) - \vu(t_{n+1})\|^4_{\vH^1}\|\nab\vu(s)\|^4_{\vL^2} \, ds\biggr] + Ck\|\ve^{n}\|^4_{\vL^2} \\\nonumber
			&=C(k+1)\biggl[  \int_{t_n}^{t_{n+1}} \|\nab(\vu(t_{n+1}) - \vu(s))\|^4_{\vL^2}\, ds \\\nonumber
			&\qquad+ C\int_{t_n}^{t_{n+1}} \, \|\vu(s) - \vu(t_{n+1})\|^4_{\vH^1}\|\nab\vu(t_{n+1})\|^4_{\vL^2} \, ds  \\\nonumber
			&\qquad+C\int_{t_n}^{t_{n+1}} \, \|\vu(s) - \vu(t_{n+1})\|^4_{\vH^1}\|\nab\vu(s)\|^4_{\vL^2} \, ds\biggr] + Ck\|\ve^{n}\|^4_{\vL^2}\\\nonumber
			&\qquad + \frac{1}{32}\Bigl(\|\ve^{n+1}\|^2_{\vL^2} - \|\ve^n\|^2_{\vL^2}\Bigr)^2.
		\end{align*}
		
		Similarly,  we can show
		\begin{align*}
			L_2 &\leq C(k+1) \int_{t_n}^{t_{n+1}} \|\ve^{n+1}\|^4_{\vL^2}\|\nab(\vu(s) - \vu(t_{n+1}))\|^4_{\vL^2} + Ck\|\ve^{n}\|^4_{\vL^2} \\\nonumber
			&\qquad+ \frac{1}{32}\bigl(\|\ve^{n+1}\|^2_{\vL^2} - \|\ve^n\|^2_{\vL^2}\bigr)^2.
		\end{align*}
		
		By the definition, 
		\begin{align*}
			L_3 &=  \frac{4C_L^2}{\nu}\int_{t_n}^{t_{n+1}}\|\ve^{n+1}\|^4_{\vL^2}\|\nab\vu(s)\|^2_{\vL^2}\, ds.
		\end{align*}
		We leave it alone for now, but will control it after a summation over $n$ is done later.
		
		To bound $L_4$, we again use the same rewriting of  $\|\ve^{n+1}\|^2_{\vL^2}= \bigl( \|\ve^{n+1}\|^2_{\vL^2}-\|\ve^n\|^2_{\vL^2} \bigr) +\|\ve^n\|^2_{\vL^2}$ 
		and appeal to the discrete Young's inequality to get
		\begin{align*}
			L_4 &= 2\biggl\|\int_{t_n}^{t_{n+1}}\bigl(\vG(\vu(s))-\vG(\vu(t_{n}))\bigr)\,dW(s)\biggr\|^2_{\vL^2}\Bigl(\|\ve^{n+1}\|^2_{\vL^2} - \|\ve^n\|^2_{\vL^2}\Bigr)\\\nonumber
			&\qquad+ 2\biggl\|\int_{t_n}^{t_{n+1}}\bigl(\vG(\vu(s))-\vG(\vu(t_{n}))\bigr)\,dW(s)\biggr\|^2_{\vL^2}\|\ve^{n}\|^2_{\vL^2} \\\nonumber
			&\leq C\biggl\|\int_{t_n}^{t_{n+1}}\bigl(\vG(\vu(s))-\vG(\vu(t_{n}))\bigr)\,dW(s)\biggr\|^4_{\vL^2} + \frac{1}{32}\Bigl(\|\ve^{n+1}\|^2_{\vL^2} - \|\ve^n\|^2_{\vL^2}\Bigr)^2\\\nonumber
			&\qquad+ 2\biggl\|\int_{t_n}^{t_{n+1}}\bigl(\vG(\vu(s))-\vG(\vu(t_{n}))\bigr)\,dW(s)\biggr\|^2_{\vL^2}\|\ve^{n}\|^2_{\vL^2}\\\nonumber
			&= C\biggl\|\int_{t_n}^{t_{n+1}}\bigl(\vG(\vu(s))-\vG(\vu(t_{n}))\bigr)\,dW(s)\biggr\|^4_{\vL^2} + \frac{1}{32}\Bigl(\|\ve^{n+1}\|^2_{\vL^2} - \|\ve^n\|^2_{\vL^2}\Bigr)^2\\\nonumber
			&\qquad+ 2\biggl\|\int_{t_n}^{t_{n+1}}\bigl(\vG(\vu(s))-\vG(\vu(t_{n}))\bigr)\,dW(s)\biggr\|^2_{\vL^2}\|\ve^{n}\|^2_{\vL^2} \\\nonumber&\qquad- 2\int_{t_n}^{t_{n+1}}\|\vG(\vu(s)) - \vG(\vu(t_{n}))\|^2_{\vL^2}\, ds\|\ve^n\|^2_{\vL^2}\\\nonumber
			&\qquad+ 2\int_{t_n}^{t_{n+1}}\|\vG(\vu(s)) - \vG(\vu(t_{n}))\|^2_{\vL^2}\, ds\|\ve^n\|^2_{\vL^2}\\\nonumber
			&\leq C\biggl\|\int_{t_n}^{t_{n+1}}\bigl(\vG(\vu(s))-\vG(\vu(t_{n}))\bigr)\,dW(s)\biggr\|^4_{\vL^2} + \frac{1}{32}\Bigl(\|\ve^{n+1}\|^2_{\vL^2} - \|\ve^n\|^2_{\vL^2}\Bigr)^2\\\nonumber
			&\qquad+ 2\biggl\|\int_{t_n}^{t_{n+1}}\bigl(\vG(\vu(s))-\vG(\vu(t_{n}))\bigr)\,dW(s)\biggr\|^2_{\vL^2}\|\ve^{n}\|^2_{\vL^2} \\\nonumber&\qquad- 2\int_{t_n}^{t_{n+1}}\|\vG(\vu(s)) - \vG(\vu(t_{n}))\|^2_{\vL^2}\, ds\|\ve^n\|^2_{\vL^2}\\\nonumber
			&\qquad+ C\int_{t_n}^{t_{n+1}}\|\vG(\vu(s)) - \vG(\vu(t_{n}))\|^4_{\vL^2}\, ds + Ck\|\ve^n\|^4_{\vL^2},
		\end{align*}
		We note that the rewriting strategy allows us to set up the stage for applying the stochastic Gronwall inequality later.
		
		Similarly, we control $L_5 + L_6$ as follows:
		\begin{align*}
			L_5 + L_6 &=2\|(\vG(\vu(t_n)) - \vG(\vu^n))\Delta W_{n+1}\|^2_{\vL^2}\Bigl(\|\ve^{n+1}\|^2_{\vL^2}-\|\ve^n\|^2_{\vL^2}\Bigr)  \\\nonumber
			&\qquad + \biggl(\int_{t_n}^{t_{n+1}}\bigl(\vG(\vu(s))-\vG(\vu^n)\bigr)\,dW(s),\ve^n\biggr)\Bigl(\|\ve^{n+1}\|^2_{\vL^2} - \|\ve^n\|^2_{\vL^2}\Bigr)\\\nonumber
			&\qquad+2\|(\vG(\vu(t_n)) - \vG(\vu^n))\Delta W_{n+1}\|^2_{\vL^2}\|\ve^{n}\|^2_{\vL^2}  \\\nonumber
			&\qquad + \biggl(\int_{t_n}^{t_{n+1}}\bigl(\vG(\vu(s))-\vG(\vu^n)\bigr)\,dW(s),\ve^n\biggr)\|\ve^{n}\|^2_{\vL^2}\\\nonumber
			&\leq \frac{1}{32} \Bigl(\|\ve^{n+1}\|^2_{\vL^2} - \|\ve^n\|^2_{\vL^2}\Bigr)^2 + C\|\vG(\vu(t_{n})) - \vG(\vu^n)\|^4_{\vL^2}|\Delta W_{n+1}|^4 \\\nonumber
			&\qquad+ C\biggl\|\int_{t_n}^{t_{n+1}}\bigl(\vG(\vu(s)) - \vG(\vu^n)\bigr)\, dW(s)\biggr\|^2_{\vL^2}\|\ve^n\|^2_{\vL^2}\\\nonumber
			&\qquad+2\|(\vG(\vu(t_n)) - \vG(\vu^n))\Delta W_{n+1}\|^2_{\vL^2}\|\ve^{n}\|^2_{\vL^2}  \\\nonumber
			&\qquad + \biggl(\int_{t_n}^{t_{n+1}}\bigl(\vG(\vu(s))-\vG(\vu^n)\bigr)\,dW(s),\ve^n\biggr)\|\ve^{n}\|^2_{\vL^2}\\\nonumber
			&= \frac{1}{32} \Bigl(\|\ve^{n+1}\|^2_{\vL^2} - \|\ve^n\|^2_{\vL^2}\Bigr)^2 + C\|\vG(\vu(t_{n})) - \vG(\vu^n)\|^4_{\vL^2}|\Delta W_{n+1}|^4 \\\nonumber
			&\qquad+ C\biggl\|\int_{t_n}^{t_{n+1}}\bigl(\vG(\vu(s)) - \vG(\vu(t_n))\bigr)\, dW(s)\biggr\|^2_{\vL^2}\|\ve^n\|^2_{\vL^2}\\\nonumber
			&\qquad+(C+2)\|(\vG(\vu(t_n)) - \vG(\vu^n))\Delta W_{n+1}\|^2_{\vL^2}\|\ve^{n}\|^2_{\vL^2}  \\\nonumber
			&\qquad + \biggl(\int_{t_n}^{t_{n+1}}\bigl(\vG(\vu(s))-\vG(\vu^n)\bigr)\,dW(s),\ve^n\biggr)\|\ve^{n}\|^2_{\vL^2}\\\nonumber
			&\leq \frac{1}{32} \Bigl(\|\ve^{n+1}\|^2_{\vL^2} - \|\ve^n\|^2_{\vL^2}\Bigr)^2 + C\|\vG(\vu(t_{n})) - \vG(\vu^n)\|^4_{\vL^2}|\Delta W_{n+1}|^4 \\\nonumber
			&\qquad+ C\biggl\|\int_{t_n}^{t_{n+1}}\bigl(\vG(\vu(s)) - \vG(\vu(t_n))\bigr)\, dW(s)\biggr\|^2_{\vL^2}\|\ve^n\|^2_{\vL^2}  + C\|\ve^{n}\|^4_{\vL^2}|\Delta W_{n+1}|^2  \\\nonumber
			&\qquad + \biggl(\int_{t_n}^{t_{n+1}}\bigl(\vG(\vu(s))-\vG(\vu^n)\bigr)\,dW(s),\ve^n\biggr)\|\ve^{n}\|^2_{\vL^2}\\\nonumber
			%
			&\leq \frac{1}{32} \Bigl(\|\ve^{n+1}\|^2_{\vL^2} - \|\ve^n\|^2_{\vL^2}\Bigr)^2 + C\|\vG(\vu(t_{n})) - \vG(\vu^n)\|^4_{\vL^2}|\Delta W_{n+1}|^4 \\\nonumber
			&\qquad+ C\biggl\|\int_{t_n}^{t_{n+1}}\Bigl(\vG(\vu(s)) - \vG(\vu(t_n))\Bigr)\, dW(s)\biggr\|^2_{\vL^2}\|\ve^n\|^2_{\vL^2}\\\nonumber
			&\qquad- C\int_{t_n}^{t_{n+1}} \|\vG(\vu(s)) - \vG(\vu(t_n))\|^2_{\vL^2}\, ds \|\ve^n\|^2_{\vL^2}\\\nonumber
			&\qquad+ C\|\ve^{n}\|^4_{\vL^2}|\Delta W_{n+1}|^2 
			+ C\int_{t_n}^{t_{n+1}} \|\vG(\vu(s)) - \vG(\vu(t_n))\|^4_{\vL^2}\, ds  { + Ck\|\ve^n\|^4_{\vL^2} }
			\\\nonumber
			&\qquad + \biggl(\int_{t_n}^{t_{n+1}}\bigl(\vG(\vu(s))-\vG(\vu^n)\bigr)\,dW(s),\ve^n\biggr)\|\ve^{n}\|^2_{\vL^2}\\\nonumber
			&\leq \frac{1}{32} \Bigl(\|\ve^{n+1}\|^2_{\vL^2} - \|\ve^n\|^2_{\vL^2}\Bigr)^2  + Ck^2\|\ve^n\|^4_{\vL^2} \\\nonumber
			&\qquad+ C\biggl\|\int_{t_n}^{t_{n+1}}\bigl(\vG(\vu(s)) - \vG(\vu(t_n))\bigr)\, dW(s)\biggr\|^2_{\vL^2}\|\ve^n\|^2_{\vL^2}\\\nonumber
			&\qquad { - C\int_{t_n}^{t_{n+1}} \|\vG(\vu(s)) - \vG(\vu(t_n))\|^2_{\vL^2}\, ds \|\ve^n\|^2_{\vL^2} } \\\nonumber
			&\qquad+ C\|\ve^{n}\|^4_{\vL^2}|\Delta W_{n+1}|^2  { + Ck\|\ve^n\|^4_{\vL^2}  }
			+ C\int_{t_n}^{t_{n+1}} \|\vG(\vu(s)) - \vG(\vu(t_n))\|^4_{\vL^2}\, ds   \\\nonumber
			&\qquad+ \biggl(\int_{t_n}^{t_{n+1}}\bigl(\vG(\vu(s))-\vG(\vu^n)\bigr)\,dW(s),\ve^n\biggr)\|\ve^{n}\|^2_{\vL^2}.
		\end{align*}
		
		\smallskip
		{\em Step 3:}  	Substituting all the estimates for $L_1, \cdots, L_6$ into the right-hand side of \eqref{eq_3.17} and combining the like-terms, we obtain  
		\begin{align}\label{eq_3.18}
			&\frac14\bigl[\|\ve^{n+1}\|^4_{\vL^2} - \|\ve^n\|^4_{\vL^2}\bigr] + \frac{3\nu}{8}\|\nab\ve^{n+1}\|^2_{\vL^2}\|\ve^{n+1}\|^2_{\vL^2} \\\nonumber
			&\qquad\quad   + \frac18\bigl(\|\ve^{n+1}\|^2_{\vL^2} - \|\ve^n\|^2_{\vL^2}\bigr)^2 
			+ \frac14\|\ve^{n+1} - \ve^n\|^2_{\vL^2}\|\ve^{n+1}\|^2_{\vL^2} \\\nonumber 
			&\quad \leq C(k+1)\biggl[  \int_{t_n}^{t_{n+1}} \|\nab(\vu(t_{n+1}) - \vu(s))\|^4_{\vL^2}\, ds \\\nonumber
			&\qquad\quad + C\int_{t_n}^{t_{n+1}} \, \|\vu(s) - \vu(t_{n+1})\|^4_{\vH^1}\|\nab\vu(t_{n+1})\|^4_{\vL^2} \, ds  \\\nonumber
			&\qquad\quad +C\int_{t_n}^{t_{n+1}} \, \|\vu(s) - \vu(t_{n+1})\|^4_{\vH^1}\|\nab\vu(s)\|^4_{\vL^2} \, ds\biggr] 
			\\\nonumber
			&\qquad\quad  +C(k+1) \int_{t_n}^{t_{n+1}} \|\ve^{n+1}\|^4_{\vL^2}\|\nab(\vu(s) - \vu(t_{n+1}))\|^4_{\vL^2} { + 2Ck\|\ve^{n}\|^4_{\vL^2} } +L_3\\\nonumber
			&\qquad\quad +C\biggl\|\int_{t_n}^{t_{n+1}}\bigl(\vG(\vu(s))-\vG(\vu(t_{n}))\bigr)\,dW(s)\biggr\|^4_{\vL^2} \\\nonumber
			&\qquad\quad 
			+2\biggl\|\int_{t_n}^{t_{n+1}}\bigl(\vG(\vu(s))-\vG(\vu(t_{n}))\bigr)\,dW(s)\biggr\|^2_{\vL^2}\|\ve^{n}\|^2_{\vL^2} \\\nonumber
			&\qquad\quad - 2\int_{t_n}^{t_{n+1}}\|\vG(\vu(s)) - \vG(\vu(t_{n}))\|^2_{\vL^2}\, ds\|\ve^n\|^2_{\vL^2}\\\nonumber
			&\qquad\quad + C\int_{t_n}^{t_{n+1}}\|\vG(\vu(s)) - \vG(\vu(t_{n}))\|^4_{\vL^2}\, ds + Ck\|\ve^n\|^4_{\vL^2}\\\nonumber
			&\qquad\quad + C\biggl\|\int_{t_n}^{t_{n+1}}\bigl(\vG(\vu(s)) - \vG(\vu(t_n))\bigr)\, dW(s)\biggr\|^2_{\vL^2}\|\ve^n\|^2_{\vL^2}\\\nonumber
			&\qquad\quad { - C\int_{t_n}^{t_{n+1}} \|\vG(\vu(s)) - \vG(\vu(t_n))\|^2_{\vL^2}\, ds \|\ve^n\|^2_{\vL^2} } \\\nonumber
			&\qquad\quad + C\|\ve^{n}\|^4_{\vL^2}|\Delta W_{n+1}|^2  
			+ C\int_{t_n}^{t_{n+1}} \|\vG(\vu(s)) - \vG(\vu(t_n))\|^4_{\vL^2}\, ds  
			\\\nonumber
			&\qquad\quad  + \biggl(\int_{t_n}^{t_{n+1}}\bigl(\vG(\vu(s))-\vG(\vu^n)\bigr)\,dW(s),\ve^n\biggr)\|\ve^{n}\|^2_{\vL^2}\\\nonumber
			&\quad \leq \biggl\{  C\int_{t_n}^{t_{n+1}} \|\nab(\vu(t_{n+1}) - \vu(s))\|^4_{\vL^2}\, ds \\\nonumber
			&\qquad\quad + C\int_{t_n}^{t_{n+1}} \, \|\vu(s) - \vu(t_{n+1})\|^4_{\vH^1}\|\nab\vu(t_{n+1})\|^4_{\vL^2} \, ds  \\\nonumber
			&\qquad\quad +C\int_{t_n}^{t_{n+1}} \, \|\vu(s) - \vu(t_{n+1})\|^4_{\vH^1}\|\nab\vu(s)\|^4_{\vL^2} \, ds  \\\nonumber
			&\qquad\quad +C\int_{t_n}^{t_{n+1}} \|\ve^{n+1}\|^4_{\vL^2}\|\nab(\vu(s) - \vu(t_{n+1}))\|^4_{\vL^2} \\\nonumber
			&\qquad\quad+C\biggl\|\int_{t_n}^{t_{n+1}}\bigl(\vG(\vu(s))-\vG(\vu(t_{n}))\bigr)\,dW(s)\biggr\|^4_{\vL^2}\\\nonumber
			&\qquad\quad +C\int_{t_n}^{t_{n+1}}\|\vG(\vu(s)) - \vG(\vu(t_{n}))\|^4_{\vL^2}\, ds 
			\biggr\}\\\nonumber
			&\qquad \quad +\biggl\{(C+2)\biggl\|\int_{t_n}^{t_{n+1}}\bigl(\vG(\vu(s))-\vG(\vu(t_{n}))\bigr)\,dW(s)\biggr\|^2_{\vL^2}\|\ve^{n}\|^2_{\vL^2} \\\nonumber
			&\qquad\quad {- (C+2)\int_{t_n}^{t_{n+1}}\|\vG(\vu(s)) - \vG(\vu(t_{n}))\|^2_{\vL^2}\, ds\|\ve^n\|^2_{\vL^2} } \\\nonumber
			&\qquad\quad + C\|\ve^{n}\|^4_{\vL^2}|\Delta W_{n+1}|^2 { - Ck\|\ve^n\|^4_{\vL^2} } \\ \nonumber
			&\qquad \quad
			+ \biggl(\int_{t_n}^{t_{n+1}}\bigl(\vG(\vu(s))-\vG(\vu^n)\bigr)\,dW(s),\ve^n\biggr)\|\ve^{n}\|^2_{\vL^2}\biggr\} \\\nonumber 
			&\qquad\quad {  + \Bigl(3Ck\|\ve^n\|^4_{\vL^2} + L_3\Bigr).}
		\end{align}
		
		Next, applying the summation operator $\sum_{n=0}^{\ell}$ to \eqref{eq_3.18} for $0<\ell \leq M-1$ and using the following estimate for the term $L_3$ 
		\begin{align*}
			\sum_{n = 0}^{\ell}L_3 
			&=  \frac{4C_L^2}{\nu}\sum_{n = 0}^{\ell-1} \|\ve^{n+1}\|^4_{\vL^2} \int_{t_{n}}^{t_{n+1}} \|\nab\vu(s)\|^2_{\vL^{2}}\, ds  +  \frac{4C_L^2}{\nu}\|\ve^{\ell}\|^4_{\vL^2}\int_{t_{\ell}}^{t_{\ell+1}}\|\nab\vu(s)\|^2_{\vL^2}\, ds \\\nonumber
			&\quad+  \frac{4C_L^2}{\nu}\Bigl[\|\ve^{\ell +1}\|^2_{\vL^2} + \|\ve^{\ell }\|^2_{\vL^2}\Bigr] \Bigl[\|\ve^{\ell +1}\|^2_{\vL^2}-\|\ve^{\ell}\|^2_{\vL^2}\Bigr] \int_{t_{\ell}}^{t_{\ell+1}}\|\nab\vu(s)\|^2_{\vL^2}\, ds  \\\nonumber
			&\leq  \frac{4C_L^2}{\nu}\sum_{n =0}^{\ell} \|\ve^{n}\|^4_{\vL^2} \int_{t_{n}}^{t_{n+1}}\|\nab\vu(s)\|^2_{\vL^2}\, ds \\\nonumber
			&+ Ck  \Bigl[\|\ve^{\ell +1}\|^4_{\vL^2} + \|\ve^{\ell }\|^4_{\vL^2}\Bigr] \int_{t_{\ell}}^{t_{\ell+1}}\|\nab\vu(s)\|^4_{\vL^2}\, ds + { \frac{1}{32}\sum_{n=0}^{\ell}\Bigl(\|\ve^{n+1}\|^2_{\vL^2} - \|\ve^{n}\|^2_{\vL^2}\Bigr)^2},
		\end{align*}
		we obtain the following master inequality:
		\begin{align}\label{eq3.19}
			&\frac14\|\ve^{\ell+1}\|^4_{\vL^2} + \frac{3}{32}\sum_{n = 0}^{\ell}\bigl(\|\ve^{n+1}\|^2_{\vL^2} - \|\ve^n\|^2_{\vL^2}\bigr)^2 + \frac{3\nu}{8}\sum_{n = 0}^{\ell}\|\nab\ve^{n+1}\|^2_{\vL^2}\|\ve^{n+1}\|^2_{\vL^2}  \\\nonumber
			&\hskip 0.99in  + \frac14\sum_{n = 0}^{\ell}\|\ve^{n+1} - \ve^n\|^2_{\vL^2}\|\ve^{n+1}\|^2_{\vL^2}  
			\leq \mathcal{D}_{\ell} + Y_{\ell} + \sum_{n = 0}^{\ell} \mathcal{S}_{n}\|\ve^n\|^4_{\vL^2},
		\end{align}
		where 
		\begin{align*}
			\mathcal{D}_{\ell} &:= C\sum_{n = 0}^{\ell} \biggl\{  \int_{t_n}^{t_{n+1}} \|\nab(\vu(t_{n+1}) - \vu(s))\|^4_{\vL^2}\, ds \\\nonumber
			&\qquad+ \|\nab\vu(t_{n+1})\|^4_{\vL^2} \int_{t_n}^{t_{n+1}} \, \|\vu(s) - \vu(t_{n+1})\|^4_{\vH^1} \, ds  \\\nonumber&\qquad+\int_{t_n}^{t_{n+1}} \, \|\vu(s) - \vu(t_{n+1})\|^4_{\vH^1}\|\nab\vu(s)\|^4_{\vL^2} \, ds  \\\nonumber
			&\qquad +\|\ve^{n+1}\|^4_{\vL^2}\int_{t_n}^{t_{n+1}} \|\nab(\vu(s) - \vu(t_{n+1}))\|^4_{\vL^2}\, ds \\\nonumber
			&\qquad + \biggl\|\int_{t_n}^{t_{n+1}}\bigl(\vG(\vu(s))-\vG(\vu(t_{n}))\bigr)\,dW(s)\biggr\|^4_{\vL^2}\\\nonumber
			&\qquad+\int_{t_n}^{t_{n+1}}\|\vG(\vu(s)) - \vG(\vu(t_{n}))\|^4_{\vL^2}\, ds 
			\biggr\} \\\nonumber
			&\qquad+Ck \Bigl[\|\ve^{\ell +1}\|^4_{\vL^2} + \|\ve^{\ell }\|^4_{\vL^2}\Bigr] \int_{t_{\ell}}^{t_{\ell+1}}\|\nab\vu(s)\|^4_{\vL^2}\, ds ,\\\nonumber
			\mathcal{S}_n &:= \frac{4C_L^2}{\nu}\int_{t_{n}}^{t_{n+1}} \|\nab \vu(s)\|^2_{\vL^2} \, ds + Ck,  \qquad 
			R_{\ell} := \sum_{n=0}^{\ell} Y_{n},
		\end{align*}
		and 
		\begin{align*}
			Y_{n} &:=(C+2)\|\ve^{n}\|^2_{\vL^2} \biggl\|\int_{t_n}^{t_{n+1}}\bigl(\vG(\vu(s))-\vG(\vu(t_{n}))\bigr)\,dW(s)\biggr\|^2_{\vL^2} \\\nonumber&\qquad- (C+2)\|\ve^n\|^2_{\vL^2} \int_{t_n}^{t_{n+1}}\|\vG(\vu(s)) - \vG(\vu(t_{n}))\|^2_{\vL^2}\, ds\\\nonumber
			&\qquad+ C\|\ve^{n}\|^4_{\vL^2}|\Delta W_{n+1}|^2 - Ck\|\ve^n\|^4_{\vL^2}  \\\nonumber
			&\qquad+\|\ve^{n}\|^2_{\vL^2}\biggl(\int_{t_n}^{t_{n+1}}\bigl(\vG(\vu(s))-\vG(\vu^n)\bigr)\,dW(s),\ve^n\biggr).			 
		\end{align*}
		
		\smallskip
		{\em Step 4:}  Suppose that $\{R_{\ell}; \ell\geq 0\}$  is a martingale  (which will be verified in {\em Step 5} below), applying the stochastic Gronwall inequality  \eqref{ineq2.4} to \eqref{eq3.19}, we get
		\begin{align}\label{eq_3.20}
			\Bigl(\mE\bigl[\sup_{1 \leq n \leq M}\|\ve^n\|^{4q}_{\vL^2}\bigr]\Bigr)^{\frac{1}{4q}} 
			&\leq \left(1+ \frac{1}{1- \alpha q}\right)^{\frac{1}{4\alpha q}}  \left\{\mE\left[\exp\left(\beta q \sum_{n = 0}^{M-1}\mathcal{S}_{n}\right)\right]\right\}^{\frac{1}{4 \beta q}}   \\\nonumber
			&\qquad\times\left\{\mE\left[\sup_{0 \leq \ell \leq M-1} \mathcal{D}_{\ell}\right]\right\}^{\frac14},
		\end{align}
		with $\alpha =1+\eps$, $\beta = \frac{1+\eps}{\eps}$, and $0<q\leq 1-\eps$.		
		
		To prove the desired error estimate, we need to control the two expectation factors on the right-hand side of \eqref{eq_3.20}. Similar to the arguments used to bound \eqref{eq_3.14},  it follows from Lemma \ref{lemma_exp} that 
		\begin{align*}
			\mE\left[\exp\left(\beta q \sum_{n = 0}^{M-1}\mathcal{S}_{n}\right)\right]
			&\leq  \mE\left[\exp\left(\beta q \bigg[CT  +  \frac{4C_L^2}{\nu}\int_{0}^{T}\|\nab\vu(s)\|^2_{\vL^2}\, ds\bigg]  \right)\right] \leq C_3e^{C\beta q T}.
		\end{align*}
		Moreover, the second factor can be bounded as follows:
		\begin{align*}
			\mE\left[\sup_{0 \leq \ell \leq M-1} \mathcal{D}_{\ell}\right] 
			&\leq  C\sum_{n = 0}^{M-1}\mE \biggl[  \int_{t_n}^{t_{n+1}} \|\nab(\vu(t_{n+1}) - \vu(s))\|^4_{\vL^2}\, ds \\\nonumber
			&\qquad+ \|\nab\vu(t_{n+1})\|^4_{\vL^2} \int_{t_n}^{t_{n+1}} \, \|\vu(s) - \vu(t_{n+1})\|^4_{\vH^1}\, ds  \\\nonumber&\qquad+\int_{t_n}^{t_{n+1}} \, \|\vu(s) - \vu(t_{n+1})\|^4_{\vH^1}\|\nab\vu(s)\|^4_{\vL^2} \, ds  \\\nonumber
			&\qquad +\|\ve^{n+1}\|^4_{\vL^2} \int_{t_n}^{t_{n+1}} |\nab(\vu(s) - \vu(t_{n+1}))\|^4_{\vL^2}\, ds \\\nonumber
			&\qquad + \biggl\|\int_{t_n}^{t_{n+1}}\bigl(\vG(\vu(s))-\vG(\vu(t_{n}))\bigr)\,dW(s)\biggr\|^4_{\vL^2}\\\nonumber
			&\qquad+\int_{t_n}^{t_{n+1}}\|\vG(\vu(s)) - \vG(\vu(t_{n}))\|^4_{\vL^2}\, ds 
			\biggr] \\\nonumber
			&\qquad+Ck^2 \mE\left[\sup_{s\in [0,T]}\|\nab\vu(s)\|^4_{\vL^2}\, \sup_{0 \leq \ell \leq M-1}\left[\|\ve^{\ell +1}\|^4_{\vL^2} + \|\ve^{\ell }\|^4_{\vL^2}\right]\right]\\\nonumber
			&=: \mathcal{U}_1 + \mathcal{U}_2.
		\end{align*}
		
		Using the same techniques as employed for the estimations of $T_1,\cdots, T_5$ in the proof of Theorem \ref{theorem_semi_chapter5}, we can obtain 
		\begin{align*}
			\mathcal{U}_1 \leq C\Bigl(C_{2,4} + C_{2,8}C_{1,8} + C_{2,8}(C_{1,8} + C_{4,3}+ C_{\vG} C_{2,4} \Bigr)k^{4\gamma} \quad \mbox{for any } 0<\gamma< \frac12.
		\end{align*}
		It follows from Schwarz inequality, Lemmas \ref{lemm3.3} and \ref{stability_means} that
		\begin{align*}
			\mathcal{U}_2 &= Ck^2 \mE\left[\sup_{s\in [0,T]}\|\nab\vu(s)\|^4_{\vL^2}\, \sup_{0 \leq \ell \leq M-1}\Bigl[\|\ve^{\ell +1}\|^4_{\vL^2} + \|\ve^{\ell }\|^4_{\vL^2}\Bigr]\right]\\\nonumber
			&\leq Ck^2 \left(\mE\left[ \left(\sup_{s\in [0,T]} \|\nab \vu(s)\|^4_{\vL^2}\right)^2\right]\right)^{\frac12}\left(\mE\left[ \left(\sup_{0\leq \ell \leq M} \|\ve^{\ell}\|^4_{\vL^2}\right)^2\right]\right)^{\frac12}\\\nonumber
			&\leq CC_{1,8}(C_{1,8} + C_{4,3})k^2.
		\end{align*}
		
		Thus,  we have proved \eqref{eq3.16}  holds with $r = 2$ after substituting the above estimates for $\mathcal{U}_1$ and $\mathcal{U}_2$ into \eqref{eq_3.20}, provided that the claim holds. 
		
		\smallskip
		{\em Step 5:}  
		We now verify our claim that $\{R_{\ell}; \ell\geq 0\}$ is a martingale and proceed similarly as we did in {\em Step 7} of the proof of Theorem \ref{theorem_semi_chapter5} using Lemmas \ref{lemm3.3} and \ref{stability_means} as well as the Burkholder-Davis-Gundy inequality \eqref{BDG}.  First, we have
		\begin{align*}
			\mE[|R_{\ell}|] \leq \sum_{n = 0}^{\ell} \mE[|Y_n|] < \infty.
		\end{align*}
		Then,  it follows from It\^o's isometry and the martingale property of It\^o's integrals that  for $n\geq 0$
		\begin{align*}
			\mE[Y_n] &=  C\mE\bigg[\left\|\int_{t_n}^{t_{n+1}} (\vG(\vu(s)) - \vG(\vu(t_n)))\, dW(s)\right\|^2_{\vL^2}\|\cQ_h\ve^{n}\|^2_{\vL^2} \\\nonumber
			&\qquad- \int_{t_n}^{t_{n+1}} \|\vG(\vu(s)) - \vG(\vu(t_n))\|^2_{\vL^2}\|\cQ_h\ve^n\|^2_{\vL^2}\, ds\bigg]\\\nonumber
			&\qquad+ C\mE\bigg[\left\|\vG(\vu(t_{n})) - \vG(\vu_h^n)\right\|^2_{\vL^2}|\Delta W_{n+1}|^2\|\cQ_h\ve^{n}\|^2_{\vL^2} \\\nonumber
			&\qquad-  k\|\vG(\vu(t_{n})) - \vG(\vu_h^n)\|^2_{\vL^2}\|\cQ_h\ve^n\|^2_{\vL^2} \bigg]\\\nonumber
			&\qquad+C\mE\bigg[ \left\|\int_{t_n}^{t_{n+1}}\left(\vG(\vu(s)) - \vG(\vu_h^n)\right)\, dW(s)\right\|^2_{\vL^2} \left\|\cQ_h\ve^{n}\right\|^2_{\vL^2} \\\nonumber
			&\qquad- \int_{t_n}^{t_{n+1}} \|\vG(\vu(s)) - \vG(\vu_h^n)\|^2_{\vL^2}\|\cQ_h\ve^n\|^2_{\vL^2}\, ds\bigg]\\\nonumber
			&\qquad +\mE\bigg[ \left(\int_{t_n}^{t_{n+1}}\left(\vG(\vu(s)) - \vG(\vu_h^n)\right)\, dW(s), \cQ_h\ve^{n}\right)\|\cQ_h\ve^{n}\|^2_{\vL^2}\bigg] =0.
		\end{align*}
		Moreover,  the conditional expectation of $R_{\ell}$ given $\{R_n\}_{n=0}^{\ell-1}$ satisfies
		\begin{align*}
			\mE[R_{\ell}| R_0, R_1, \cdots, R_n] 
			&= \mE[Y_0+Y_1 + \cdots + Y_n| R_0, R_1, \cdots, R_n] \\ 
			&\qquad + \mE[Y_{n+1}  + \cdots+Y_{\ell}| R_0, R_1,\cdots, R_{n}]\\\nonumber
			&= R_{n} + \mE[Y_{n+1} + \cdots + Y_{\ell}] = R_n.
		\end{align*}
		Thus, $\bigl\{R_{\ell}\bigr\}_{\ell \geq 0}$ is a martingale. 
		
		\smallskip
		{\em Step 6:}	To prove \eqref{eq3.16} for $r>2$, we use an inductive argument as it was done in the proof of \cite[Theorem 3.3 ]{LV2021}, starting with multiplying the error inequality \eqref{eq_3.12} by $\|\ve^{n+1}\|^{2^j}_{\vL^2}$ (assuming \eqref{eq3.16} holds for $r=j$),
		and utilize the technique of leveraging the stochastic Gronwall's inequality as shown for the case $r=2$. We skip the derivation to save space and leave it to the interested reader to complete. 
	\end{proof}
	
	In the corollary below, we present a high moment error estimate of Theorem \ref{theorem_semi_higher_moment}  in a familiar form and also state a pathwise error estimate in the  $L^2$-norm. 
	
	\begin{corollary}\label{theorem_velocity} 
		Let $0<\gamma<\frac12$ be the same as in Theorem \ref{theorem_semi_higher_moment}. Then there hold the following estimates:  
		\begin{subequations}\label{eq_3.23}
			\begin{align}\label{eq_3.23a}
				&\left(\mE\left[\max_{1\leq n\leq M}\|\vu(t_n) - \vu^n\|^{m}_{\vL^2}\right]\right)^{\frac{1}{m}} 
				\leq \tilde{C}_2\, k^{\gamma} \qquad\qquad\forall m>2,\\
				&\max_{1\leq n\leq M}\|\vu(t_n) - \vu^n\|_{\vL^2}
				\leq K(\omega)\, k^{\gamma_1} \hskip 0.95in\mbox{$\mathbb{P}$-a.s.,} \label{eq_3.23b}
			\end{align}
		\end{subequations}
		where  $\tilde{C}_2$ is a positive constant,   $K$ is a random variable such that $E[|K|^m]<\infty$, and $0<\gamma_1<\gamma-\frac{1}{m}$ provided that $m> \frac{1}{\gamma}$. 
	\end{corollary}
	
	\begin{proof}
		To show \eqref{eq_3.23a},  taking $q = \frac12$ in \eqref{eq3.16} yields
		\begin{align*}
			\Bigl(\mE\bigl[\max_{1\leq n\leq M}\|\vu(t_n) - \vu^n\|^{2^{r-1}}_{\vL^2}\bigr]\Bigr)^{\frac{1}{2^{r-1}}} 
			\leq C\, k^{\gamma}\qquad\qquad \forall r\in \mathbb{N}.
		\end{align*}
		\eqref{eq_3.23a} follows from the above estimate and an application of the H\"older's inequality.
		
		The assertion \eqref{eq_3.23b} is an immediate consequence of  \eqref{eq_3.23a} and Kolmogorov Criteria/Theorem for a pathwise continuity of stochastic processes (cf. \cite{PZ1992}). 
	\end{proof}

	\subsection{Second moment  error estimates in the energy norm for the velocity and pressure approximations}
	The previous subsection focused on the error estimates for the velocity approximation in the $L^2$-norm; the goal of this subsection is to establish second moment error estimates for the velocity approximation in the energy norm and for the pressure approximation in a time-averaged $L^2$-norm. 
	
	\begin{theorem}\label{theorem_semi_energy} 
		Let $\vu$ be the variational solution to \eqref{equu2.100}, $\{\vu^{n}\}_{n=1}^M$ be generated by  Algorithm 1, and $\sigma_0 = \frac{1}{16K^2}$. Assume $\vu_0 \in L^{8}\bigl(\Ome;\mV\cap \vH^2(D) \bigr)\cap L^{40}\bigl(\Ome,\mV\bigr)$ such that $\mE\left[\exp\left(4\sigma\|\vu_0\|^2_{\vL^2}\right)\right] < \infty$ for any $\sigma \in (0,\sigma_0]$, and $\vG$ satisfies Assumptions (B1)--(B3). Then, for any $0 < \gamma < \frac12$, there holds
		\begin{align}\label{eqn_5.21}
			\biggl(\mE\biggl[\nu k \sum_{n = 1}^M \|\nab(\vu(t_n) - \vu^n)\|^2_{\vL^2}\biggr]\biggr)^{\frac{1}{2}}  			
			\leq \tilde{C}_3\, k^{\gamma}
		\end{align}
		for some positive constant $\tilde{C}_3 = C(\vu_0, T,\tilde{C}_2)$.
	\end{theorem}
	
		The above desired inequality is a direct consequence of Theorem \ref{theorem_semi_higher_moment}, Lemma \ref{lemma2.2}, and inequality \eqref{eq_3.13}. We omit the details to save space and leave it to the interested reader to verify.
	
	We conclude this section and paper by presenting the final main result of this paper, that is, to derive an error estimate for the pressure approximation in a time-averaged $L^2$-norm.  This is achieved by using the above estimate for the velocity approximation and the stochastic inf-sup condition (cf. \cite{FL2020,FPL2021}).
	
	\begin{theorem}\label{theorem_pressure} 
		Under the assumptions of Theorem \ref{theorem_semi_energy}, there holds the following estimate for the pressure approximation: 
		\begin{align}\label{eqn_5.23}
			\left(\mE\left[\left\|P(t_{\ell}) - k\sum_{n=1}^{\ell}p^{n}\right\|^2_{L^2}\right]\right)^{\frac12} \leq \tilde{C}_4k^{\gamma} \qquad\mbox{for any } 1\leq \ell \leq M,
		\end{align}
		where $\tilde{C}_4 = C(\vu_0,T,\tilde{C}_2,\tilde{C}_3,\beta)$ is a positive constant.
	\end{theorem}
	
	\begin{proof}
		The proof utilizes the stochastic version of the following inf-sup condition:
		\begin{align}\label{inf-sup}
			\sup_{\pphi \in \mH^1_{per}(D)\atop  \pphi \neq 0}\frac{\left(q, \div \pphi\right) }{\|\nab\pphi\|_{\vL^2}} \geq \beta \|q\|_{L^2}\qquad \forall q\in L^2_{per}(D)
		\end{align}
		
		Let $\E^{\ell}_P: = P(t_{\ell}) - k\sum_{n = 1}^{\ell}p^n$.  Applying the summation operator $\sum_{n=0}^\ell$ to \eqref{equu3.1} for $1\leq \ell <M$ yields 
		\begin{align}\label{eq_3.24}
			&\bigl(\vu^{\ell+1} - \vu^0,\pphi\bigr) + \nu\left(k\sum_{n = 0}^{\ell}\nab\vu^{n+1},\nab\pphi\right) +\left(k\sum_{n = 0}^{\ell}\vu^{n+1}\cdot\nab\vu^{n+1},\pphi\right) \\\nonumber
			&\hskip 0.35in	- \left(k\sum_{n = 0}^{\ell} p^{n+1},\div \pphi \right) = \left(\sum_{n = 0}^{\ell} \vG(\vu^n)\Delta W_{n+1},\pphi\right) \qquad\forall\pphi\in \vH^1_{per}(D).
		\end{align}
		Subtracting \eqref{eq_3.24} from \eqref{equu2.10a} we obtain
		\begin{align}\label{eqn_5.25}
			&\left(\ve^{\ell+1},\pphi\right) + \nu\left(k\sum_{n = 0}^{\ell} \nab\ve^{n+1},\nab\pphi\right) -\left(\E_P^{\ell+1},\div \pphi\right)\\\nonumber
			&\qquad = \nu\left(\sum_{n = 0}^{\ell}\int_{t_n}^{t_{n+1}}\nab(\vu(t_{n+1}) - \vu(s))\, ds, \nab\pphi\right) \\\nonumber
			&\qquad\qquad + \left(\sum_{n = 0}^{\ell}\int_{t_n}^{t_{n+1}}\bigl[\vu^{n+1}\cdot\nab\vu^{n+1} - \vu(t_{n+1})\cdot\nab\vu(t_{n+1})\bigr]\, ds,\pphi\right)  \\\nonumber
			&\qquad\qquad + \left(\sum_{n = 0}^{\ell}\int_{t_n}^{t_{n+1}}\bigl[\vu(t_{n+1})\cdot\nab\vu(t_{n+1}) - \vu(s)\cdot\nab\vu(s)\bigr]\, ds,\pphi\right) \\\nonumber
			&\qquad\qquad + \left(\sum_{n =0}^{\ell} \int_{t_n}^{t_{n+1}} \bigl[\vG(\vu(s)) - \vG(\vu(t_n))\bigr]\, dW(s,\pphi)\right) \\\nonumber
			&\qquad\qquad + \left(\sum_{n =0}^{\ell} \int_{t_n}^{t_{n+1}} \bigl[\vG(\vu(t_n)) - \vG(\vu^n)\bigr]\, dW(s,\pphi)\right)\\\nonumber
			&\qquad =: X_1 + X_2 + X_3 + X_4 + X_5.
		\end{align}
		
		It follows \eqref{inf-sup}, \eqref{eqn_5.25}, and Schwarz inequality  that 
		\begin{align*}
			\beta \left\| \E_P^{\ell+1}\right\|_{L^2} 
			&\leq C \|\ve^{\ell+1}\|_{\vL^2} + \nu \left\|k\sum_{n = 0}^{\ell} \nab \ve^{n+1}\right\|_{\vL^2} \\\nonumber
			&\qquad+ \frac{1}{\|\nab\pphi\|_{\vL^2}}\bigl(X_1 + X_2 + X_3 + X_4 + X_5\bigr).
		\end{align*}
		Taking expectation and using Theorem \ref{theorem_semi_higher_moment}  we get
		\begin{align*}\label{eqn_5.26}
			\beta	\mE\bigl[ \left\| \E_P^{\ell+1}\right\|^2_{L^2}\bigl]
			&\leq C k^{2\gamma}+ \mE\left[\frac{1}{\|\nab\pphi\|^2_{\vL^2}} \bigl(X_1^2+ X_2^2 + X_3^2 +X_4^2 +X_5^2\bigr) \right]
		\end{align*}
		
		Using Schwarz inequality, Theorem \ref{theorem_semi_higher_moment}, and the It\^o isometry, we  can show 
		\begin{align*}
			\mE\left[\frac{1}{\|\nab\pphi\|^2_{\vL^2}}\left\{X_1^2+ X_2^2 + X_3^2 +X_4^2 +X_5^2\right\}\right] \leq Ck^{2\gamma}.
		\end{align*}
		The proof is complete after combining the above two inequalities. 
	\end{proof}

	\bibliographystyle{abbrv}

\end{document}